\documentclass[usenatbib]{article}
\pdfoutput=1
\usepackage{graphicx}
\usepackage{amssymb,amsmath}
\usepackage{color}

\usepackage[colorinlistoftodos,shadow]{todonotes}
\usepackage{amsmath,amscd,amssymb}
\usepackage[T1]{fontenc}
\usepackage[utf8]{inputenc}
\usepackage[english]{babel}
\usepackage[numbers]{natbib}
\usepackage{color}
\usepackage{graphicx}
\usepackage{psfrag,epsfig}
\usepackage{multirow}
\usepackage{rotating}
\usepackage[letterpaper,margin=1.0in]{geometry}
\usepackage{subfigure}
\usepackage{enumerate}
\usepackage{comment}
\usepackage{lineno}
\usepackage{algorithm}
\usepackage{algpseudocode}
\usepackage{tabularx}
\usepackage[small]{caption}
\usepackage[debug=false, colorlinks=true, pdfstartview=FitV, 
linkcolor=blue, citecolor=blue, urlcolor=blue]{hyperref}
\usepackage{float}

\def\EQ#1\EN{\begin{equation}#1\end{equation}}
\def\SEQ#1\SEN{\begin{subequations}#1\end{subequations}}
\def\BA#1\EA{\begin{align}#1\end{align}}
\def\LM#1\EM{\begin{linenomath}#1\end{linenomath}}

\newcommand{\pfx}[2]{\dfrac{\partial{#1}}{\partial{#2}}}

\newcommand{\bss}[1]{\boldsymbol{#1}}

\numberwithin{equation}{section}

\title{Towards Solving the Navier-Stokes Equation on Quantum Computers}

\author{N. Ray,$^{1}$\footnote{E-mail:nray@lanl.gov}, T. Banerjee$^{2}$, B. Nadiga$^{3}$,S. Karra$^{2}$\\ 
$^{1}$Computer, Computational and Statistical Sciences (CCS-7)\\
$^{2}$Computational Earth Science Group (EES-16)\\
$^{3}$Computational Physics and Methods (CCS-2)\\
Los Alamos National Laboratory, Los Alamos, NM 87545}
\begin{document}

\maketitle

\section*{abstract}
In this paper, we explore the suitability of upcoming novel computing technologies, in particular adiabatic annealing based quantum computers, to solve fluid dynamics problems that form a critical component of several science and engineering applications. We start with simple flows with well-studied flow properties, and provide a framework to convert such systems
to a form amenable for deployment on such quantum annealers. We analyze the solutions obtained both qualitatively and quantitatively as well as the sensitivities of the various solution selection schemes on the obtained solution. \\

{\sl \noindent keywords: 
Quantum computing -- quantum annealers -- fluid dynamics -- turbulence -- linear systems.}
%
% ===============================
%
\section{Introduction}
Fluids are ubiquitous in nature and studying their dynamical properties form the core research focus for many applications. In particular, turbulent transport in many fluid flow systems underlies numerous applications such as atmospheric and climate dynamics, Inertial Confinement Fusion (ICF), combustion hydrodynamics, etc., that are of interest to academia, industry and research laboratories. Current methods for solving such systems involve numerically approximating the governing equations of flows, i.e., the Navier-Stokes (NS) equation. For example, Direct Numerical Simulations (DNS) resolve the large range of scales present in an application. Another approach called the large eddy simulations (LES) parameterizes smaller scales in the application using sub-grid scale models. The Reynolds-averaged Navier-Stokes equations (RANS) are used to describe turbulent flows based on the knowledge of the properties of turbulent statistics to approximate time-averaged solutions to the NS equations. 

Solving such complex systems on large-scale distributed machines have been extensively studied over the past decades. The computational scaling limits due to the Moore’s law on processor architecture constraints the performance of such methods on large-scale systems. In particular, 
such large-scale applications are increasingly becoming latency bound in comparison to increasing  computing power due to memory access speed and inter-processor communication latency. 
Novel computing technologies such as quantum systems are currently being explored as they provide new approaches of computation. Quantum computing theoretically promises exponential speed-up in terms of the number of states that can be explored at a time with significantly less energy requirements. 

In this paper, we explore the capability of the annealing-based quantum systems to solve fluid dynamics problems starting with transient channel flow. The DWave \citep{DWaveWS} machines were the first commercially available annealing-based quantum computers. The quantum annealer uses quantum mechanical phenomena of superposition, entanglement and tunneling to explore a given energy landscape and return a distribution of possible energy states with the global minimum ideally as the state with highest probability. 

Because of DWave's Quantum Processor Unit's (QPU) annealing based nature, it is suitable for problems of the type that minimize an unconstrained binary objective function, a.k.a., quadratic unconstrained binary optimization (QUBO) problem. The very first quantum algorithms for the DWave machines targeted optimization problems that naturally fit into the QUBO formulation. Examples include traffic route analysis \citep{10.3389/fict.2017.00029}, quantum chemistry \citep{wang2016}, etc. NP-hard problems such as finding maximum clique \citep{Chapuis2019} in graphs or graph-decomposition \citep{Ushijima-Mwesigwa:2017:GPU:3149526.3149531, MCDWave13} for which no classical polynomial time algorithms are known to exist were also studied.

Our focus is to study how a fluid dynamical system, starting with simple transient channel flows can be transformed to a form suitable for DWave's QPU, and analyse the solutions obtained when compared to classical solutions obtained using standard numerical methods. Towards that end, we use two key steps: 1) transform the problem with real data types to one with binary variables via fixed point arithmetic, and 2) pose the transformed problem as a least squares problem to convert it to a QUBO form. We finally use multiple strategies to select solutions from the distribution of the states obtained from DWave, and analyse their quantitative and qualitative properties in comparison to solutions obtained using double precision arithmetic. 

In Section \ref{background},  we begin with a brief background of governing equations of Channel flow and its numerical discretization. In Section \ref{methodology}, we describe the two key steps involved in transforming the real problem to a QUBO form along with the pre- and post processing steps. In Section \ref{nexperiments}, we discuss the numerical solutions obtained from the annealer and their analysis. Finally, we conclude with key observations in Section \ref{conclusions}.

\section{Background} \label{background}
In this section, we provide a brief overview of the discretization methodology used for numerical solution of one-dimensional time-dependent channel flow. The channel flow is a standard flow problem that is  frequently encountered in fluid mechanics. We have chosen this particular flow
as our test problem as its flow has been extensively studied and provides an excellent control for the quantum solutions. 

\subsection{Governing Equations of 1D Channel Flow}
\begin{figure}
\begin{center}
\fbox{\includegraphics[width=0.5\linewidth]%
{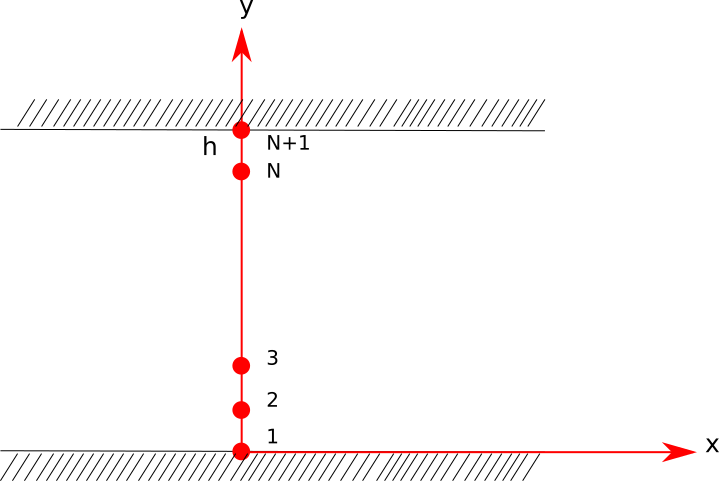}}\\[1.0ex]
\caption{1D channel flow with $N+1$ grid points.}     
\label{fig:cfschematic}    
\end{center}
\end{figure}

For the case of one-dimensional flow in a channel as shown in Figure \ref{fig:cfschematic}, that is chosen as the $x$-direction, the balance of momentum or the Navier-Stokes equation reduces to
\BA
\pfx{u_x}{t} = -\frac{1}{\rho} \pfx{\rho}{x} + \nu \pfx{^2 u_x}{y^2}, 
\label{eq:GE_1D_NS}
\EA
where $\rho$ is the fluid density, $\nu$ is the kinematic viscosity. 
Here, we assume that the $x$-direction velocity $u_x$ is only a function of the $y$-direction and that the velocity in the $y$-direction is zero. This automatically satisfies the balance of mass for incompressible flow.

We are also assuming that there are no body forces in the $x$-direction. The pressure gradient in the $x$-direction is assumed to be constant and is prescribed. In such a case, two boundary conditions (one at each channel boundary) and an initial condition for $u_x$ are needed.
We set $u_x$ to be zero at these two boundaries i.e., 
\BA
u_x(y=0) = u_x(u=h) = 0. \label{eq:BC_1D_NS}
\EA
	
\subsection{Numerical Discretization}
Finite difference methods are one of the standard numerical methods for solving PDEs. In this work, we use the second-order central difference scheme for space variables to discretize the PDE \eqref{eq:GE_1D_NS} over the channel width.  Backward/Implicit Euler scheme is used for the time evolution of the channel flow.
 
The channel is divided into $N$ intervals over $0\leq y \leq h$ with $N+1$ grid points. $u(x_i,t)$ represents the solution at the  $i$ the grid point at time  $t$ and is shortened to $u_i(t)$. Substituting the approximated operators in equations \eqref{eq:GE_1D_NS} and \eqref{eq:BC_1D_NS} leads to: 
 
\SEQ
\EQ
\frac{u_i(t+\Delta t) - u_i(t)}{\Delta t} = -\frac{1}{\rho} \pfx{\rho}{x}(t + \Delta t) 
+ \nu \frac{u_{i+1}(t+\Delta t) - 2u_i(t+\Delta t) + u_{i-1}(t+\Delta t)}{(\Delta y)^2},
\EN
\EQ
u_1(t+\Delta t) = u_{N+1}(t+\Delta t) = 0,
\EN
\SEN

Combining the equations for all the grid points, we obtain a system of linear equations: 
\EQ
\bss{A} \bss{u} = \bss{b}
\label{eq:linsysreal}
\EN

Solving \ref{eq:linsysreal} provides an approximate solution profile at time $t+\Delta t$.

\section{Solving Channel Flow on DWave Quantum Annealer} \label{methodology}
In this section we describe our methodology to convert a fluid flow problem to a form amenable for the DWave machine. An adiabatic quantum annealer can solve problems posed as unconstrained binary optimization problem. This essentially means converting the system of equations in real variables into an optimization problem with binary variables. 

Our approach uses two steps to convert the problem to QUBO. We first convert real variables to binary variables by using a fixed-point approximation. In the second step, we use a linear least-squares formulation to transform the problem to an optimization problem with binary variables. The following subsections discuss these two steps along with pre- and post-processing steps required by the annealer. 

\subsection{Fixed-point Arithmetic based Conversion of Decimal to Binary Variables}
In computing, fixed-point representation is one of the discrete representations for a real data type that has a fixed number of digits after a radix/fixed point. Essentially, it represents a real value by scaling an integer value with an implicit scaling factor which remains fixed throughout the computation. In \citep{7761616}, this idea was used to convert real variables to binary variables by using a scaling factor of 2. Thus, the solution at the \emph{i}th grid point, $u_i$ can be represented as 
\EQ
u_i = \sum_{j=1}^{n} 2^{j_0 - j}q_{j}
\label{eq:dec2bi}
\EN

Here $n$ is the precision of the representation, $j_{0}$ is the position of the fixed point, and $q_{j}$ is the $jth$ binary variable.
Clearly, keeping the precision fixed while varying $j_{0}$ leads to representing different ranges of the decimal values. Table \ref{table:prerange} shows the range of real values for various precisions with fixed point after position 1 (the leftmost position is the starting bit). While increasing the precision(\emph{n}) by one increases the range by an amount  $1/2^{n+1}$, moving the fixed point position by one would lead to doubling the maximum value. Using \ref{eq:dec2bi} to represent the solution at each grid point by using $n$ binary variables and substituting it in \ref{eq:linsysreal} leads to an extended matrix $A^{d}$ such that $A\textbf{u} = A^{d}\textbf{q}$. Finally, the linear system  \ref{eq:linsysreal} becomes 
\EQ
\bss{A^{d}} \bss{q} = \bss{b}
\label{eq:linsysbv}
\EN
Note that the right hand side of the linear system is unaffected by this transformation, and is a real data type. 

\begin{table}[]
\centering
\begin{tabular}{|c|c|c|c|c|c|c|c|c|}
 \hline
 Precision $n$ & 1 & 2 & 3 & 4 & 5 & 6 & 7 & 8  \\ \hline
 Max Real Value & 1  &  1.5 & 1.75 & 1.875 & 1.9375 & 1.96875 &  1.984375 & 1.9921875 \\ 
 \hline
\end{tabular}
 \caption{Maximum bound for reals for various precisions with fixed point after position 1.}
  \label{table:prerange}
\end{table}

\subsection{Transformation to QUBO Form}
Posing \ref{eq:linsysbv} in a least-squares form leads to the following formulation: 
\EQ
\tilde{\textbf{q}}=\underset{\textbf{q}}{min}\left\Vert A^{d}\textbf{q}-\textbf{b}\right\Vert ^{2}
\label{bilsq}
\EN
where $A^{d}$ is a real-valued matrix, \textbf{b} is a real-valued vector and \textbf{q} is a binary vector. The form of the objective function, known as the QUBO form, that the annealer can take as input is the following: 
\EQ
f\left(\textbf{q}\right)=\underset{i}{\sum}v_{i}q_{i}+\underset{i<j}{\sum}w_{ij}q_{i}q_{j}.
\EN
Here $v_{i}$ and $w_{ij}$ correspond to the weights associated with each logical qubit and the coupling strengths between two logical qubits of the problem that defines its energy landscape. In order to convert \ref{bilsq} to the QUBO form, we expand the square and use the idempotent property of binary variables to obtain
\begin{alignat*}{1}
v_{j} & =\underset{i}{\sum}A_{ij}^{d}\left(A_{ij}^{d}-2b_{i}\right)\\
w_{jk} & =2\underset{i}{\sum}A_{ij}^{d}A_{ik}^{d}
\end{alignat*}

\subsection{Pre- and Post Processing}
\begin{figure}
\centering
\subfigure[DWave Hardware Chimera Graph]{{ \includegraphics[width=.35\linewidth]{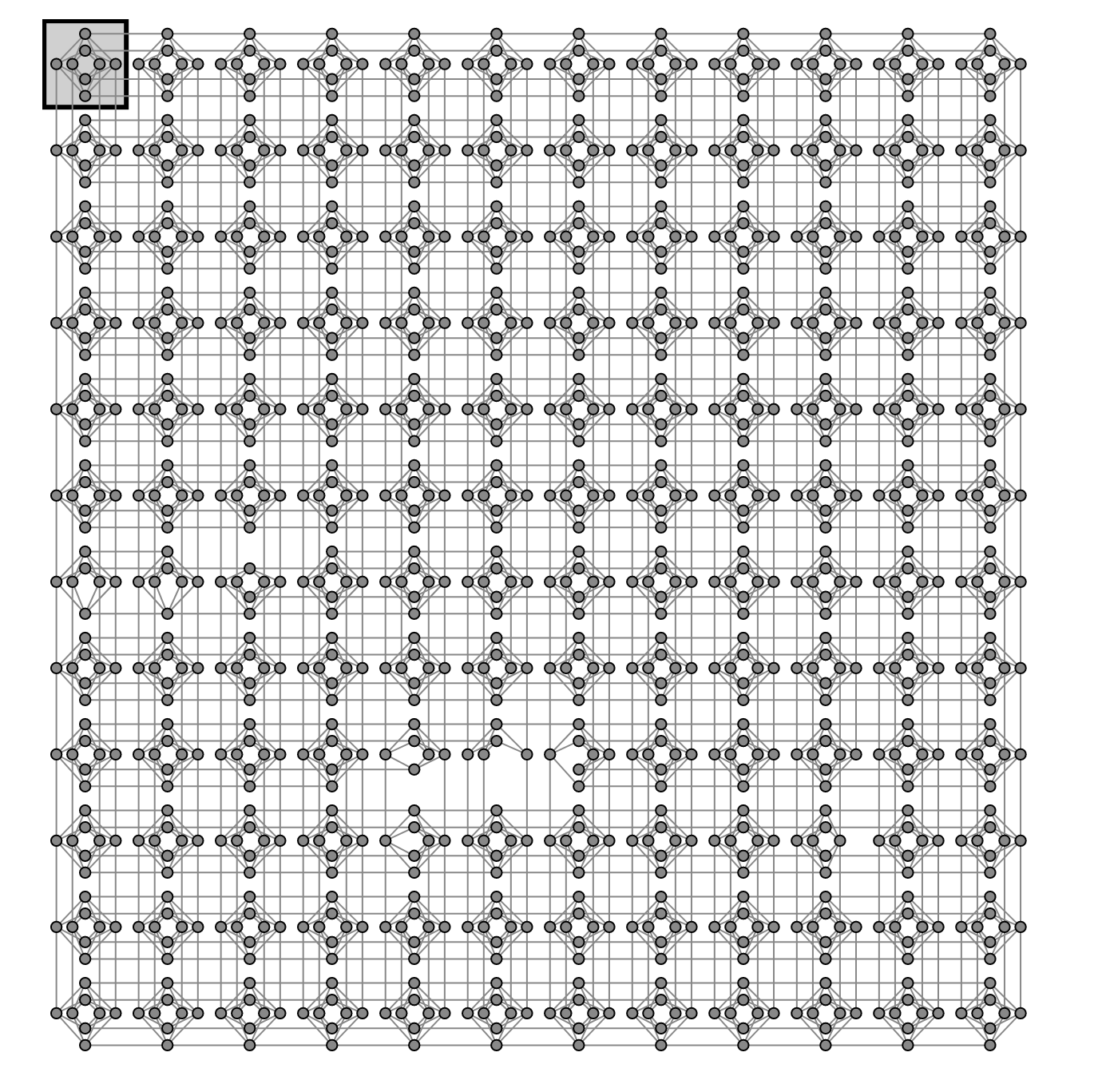} }}%
\qquad 
\subfigure[Schematic of a unit block]{{\includegraphics[width=.2\linewidth]{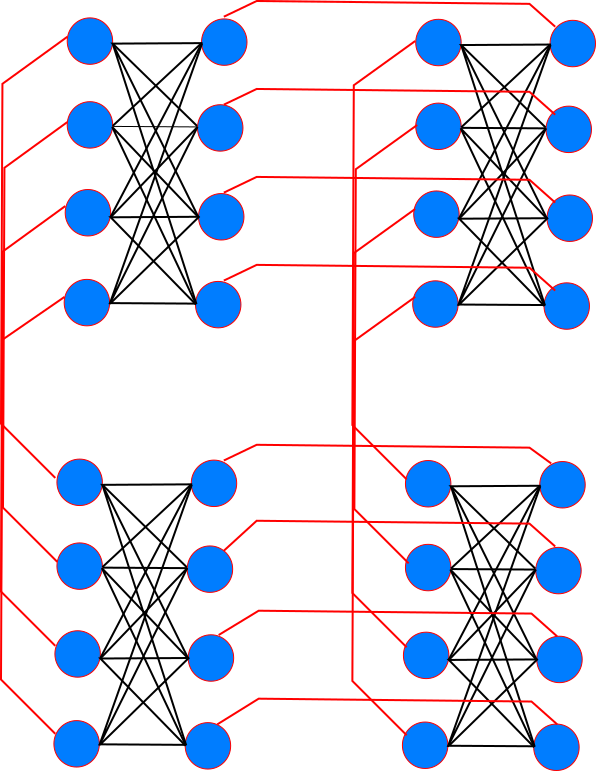} }}%
\caption{The left figure shows the layout of the physical qubits in DWave annealer. Each block in the Chimera graph is $K_{4}$ graph. The right figure shows a schematic of a 4 block configuration where the blue represents a qubit with the coupling between adjacent blocks shown in red. }
\label{fig:chimera}
\end{figure}

By design, the DWave annealer is organized as a lattice of unit blocks of qubits, where each block has eight qubits configured as a four-node bipartite graph. This hardware configuration is called a Chimera graph and is shown in Figure. \ref{fig:chimera} . As a result of this design, no three qubits are mutually coupled. This restricts the structure of the problem that the annealer can solve, and in general the logical problem needs to mapped or embedded to the hardware layout in a way that allows mutually coupled logical qubits. 

One way that DWave supports embedding generic layouts to its Chimera layout is by using the concept of \textit{chaining}. Chaining allows linking or representing a single logical qubit with a chain of hardware qubits. This involves finding sub-graphs in the Chimera layout to embed the logical problem. Such a process increases the number of qubits required for the logical problem, and as a result restricts the size of the logical problem that can be solved as we will see in the results section.

In this work, we do not focus on the problem of obtaining optimized embeddings for our channel flow problem, and instead use the utilities provided by the SAPI libraries from DWave. The Solver API (SAPI) is an interface to the DWave QPUs along with a variety of other advanced software solvers. The client applications can use this API to develop their applications in C, MATLAB, and Python. The native embedding utility provides a mapping between the logical qubits to chained hardware qubits where the algorithm tries to embed larger strongly connected components first, then smaller components. 

\begin{figure}
\begin{center}
\fbox{\includegraphics[width=1.0\linewidth]%
{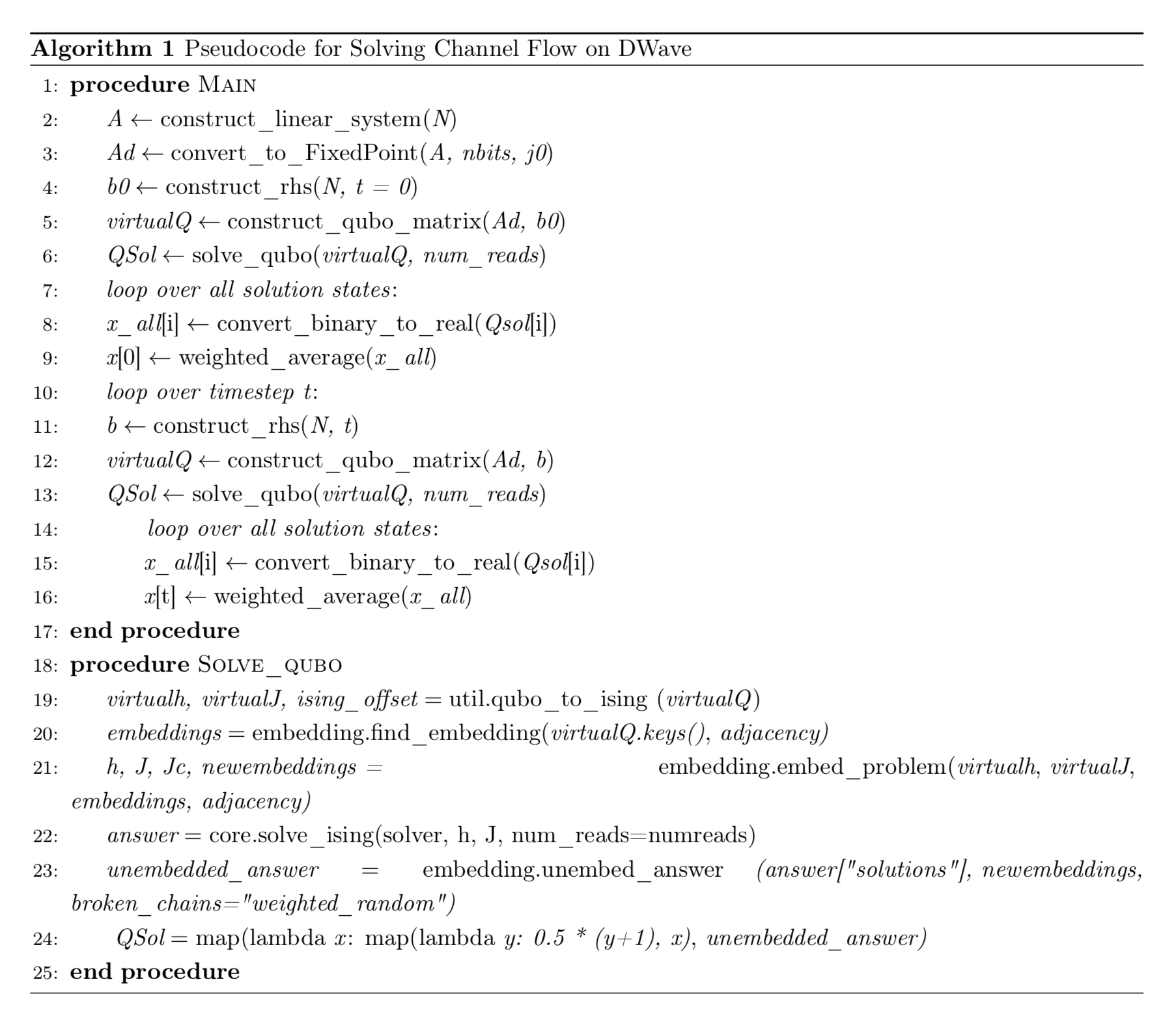}}\\[1.0ex]
\caption{Pseudocode for Solving Channel Flow on DWave}
\label{pseudocode}
\end{center}
\end{figure}

Algorithm (figure \ref{pseudocode}) provides an overview of all the steps for solving the one-dimensional channel flow on the DWave annealer. For each time step, we use the QUBO solution obtained from the previous step to construct the right-hand-side input to the linear system.  After each QUBO solve, all the states returned by the annealer are collected along with the number of occurrences of each state. Since, QUBO returns all possible solution states achieved by the hardware, we employ various strategies to compute possible solutions from the entire spectrum of returned solutions. Once a particular strategy is chosen, we transform the binary solution to a real type by using the transformation \ref{eq:dec2bi} .

In the next section, we report the results for the three key solution strategies (lowest energy, mean and weighted mean) that we use to study the effect of precision of the fixed-point representation when compared to the standard double-precision floating point arithmetic solutions for the channel flow.

%========

\section{Numerical Experiments} \label{nexperiments}
The flow parameters of the problem are shown in Table \ref{table:cfpara}. We have used the DWave DW2X\_3 hardware solver for obtaining the QUBO solutions for each time step of the flow.  The number of grid points in the mesh was varied from 5 to 10 points including the boundary points, and the simulation was performed for 10 time steps. Specifically, the QUBO solution from the previous time-step was used to form the right hand side (rhs) of the new linear system. Table. \ref{table:lpsize} shows the size of the logical QUBO with respect to the highest precision within the range 1 to 8 for which atleast one embedding was found by the DWave SAPI embedding utilities. For code development, a LANL Institutional Computing testbed, \textit{Darwin} was used along with DWave SAPI utilities (sapi-c/3.0, sapi-python/3.0, sapi-matlab/3.0, qop/2.5.0.1, etc.). In the subsequent sections, we provide both qualitative and quantitative analysis of the solutions obtain from the hardware. In particular, we explore the following three key solution selection schemes: 
\begin{itemize}
\item \textit{Lowest energy} i.e., the state with the lowest energy,  
\item \textit{Unweighted mean}, the mean of all states returned by the hardware, and 
\item \textit{Weighted mean}, the mean of all states returned by the hardware weighted by their number of occurrences.
\end{itemize}

\begin{table}
\centering
\begin{tabular}{|c|c|}
 \hline 
 Parameters:   & \\ \hline
 Channel Bounds & $[0, 1]$   \\
 Density $\rho$	 & 0.5 \\
 Viscosity $\mu$ & 0.6 \\
 Body Force $g$ & 0.4 \\
 Pressure Gradient $\delta{p}/\delta{x}$  & -2.0 \\
 Alpha $\alpha=\nu\Delta{t}/\Delta{y^2}$ & 0.4 \\
 Number of time steps & 10 \\
 \hline 
\end{tabular}
\caption{Flow parameters used for the channel flow.}
\label{table:cfpara}
\end{table}
           
\begin{table}
\centering
\begin{tabular}{|c|c|c|c|c|c|c|}
\hline
N & 5 & 6 & 7 & 8 & 9 & 10 \\ \hline
nbit & 8 & 8 & 8 & 7 & 6 & 5 \\ \hline
Logical Problem Size & 24 & 32 & 40 & 42 & 42 & 40  \\
 \hline  
\end{tabular}
\caption{Size of the logical problem with increasing number of grid points $N$ (here $N$ includes the boundary points) and precision $n$.}
\label{table:lpsize}
\end{table}

\subsection{Effect of Precision }
\begin{figure}
\begin{tabular}{|>{\raggedright}m{0.15\paperwidth}|ccc|}
\hline 
$\#$ Grid Points & 5  & 7  & 9  \tabularnewline
\cline{1-1} 
$\#$ Precision &  &  &   \tabularnewline
\hline 
2 &  {\footnotesize{}\includegraphics[scale=0.25]{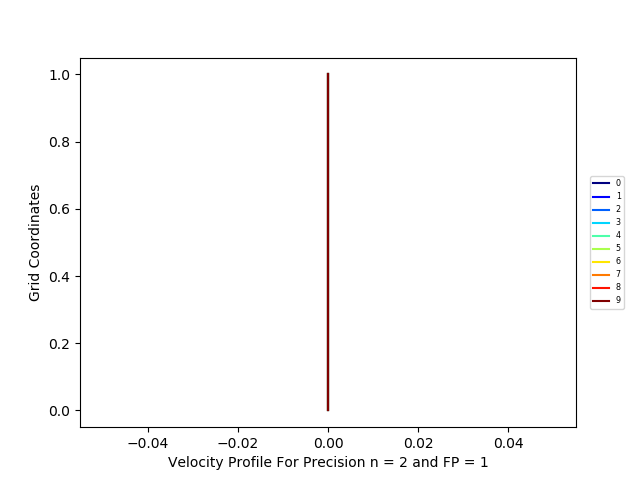}}
 & {\footnotesize{}\includegraphics[scale=0.25]{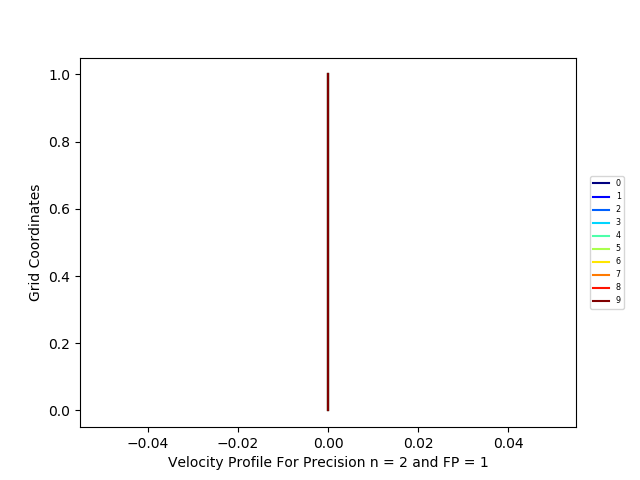}}
 & {\footnotesize{}\includegraphics[scale=0.25]{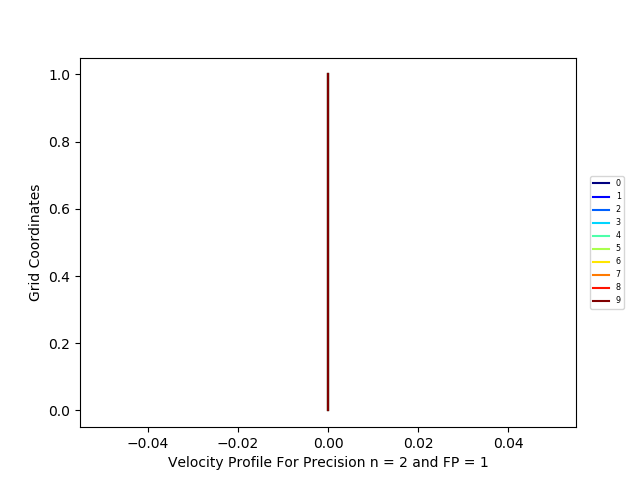}} \tabularnewline
\hline 
4 & {\footnotesize{}\includegraphics[scale=0.25]{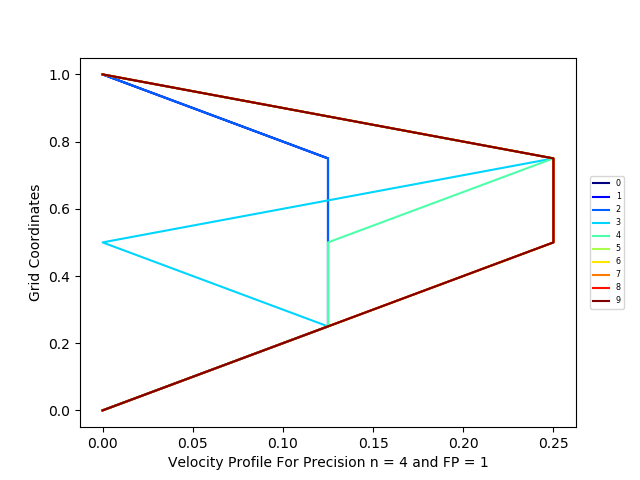}}
 & {\footnotesize{}\includegraphics[scale=0.25]{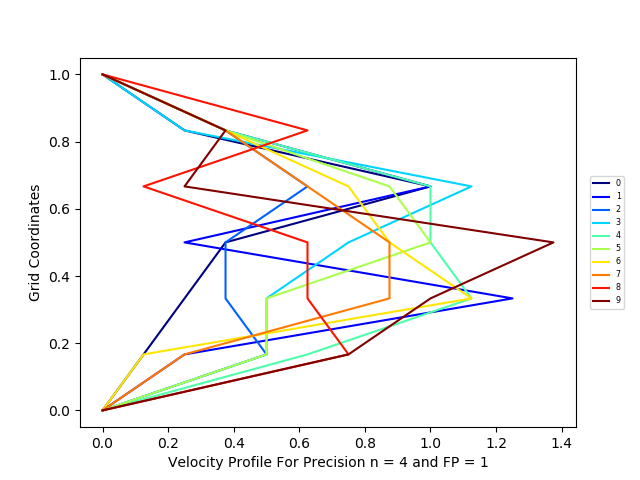}}
 & {\footnotesize{}\includegraphics[scale=0.25]{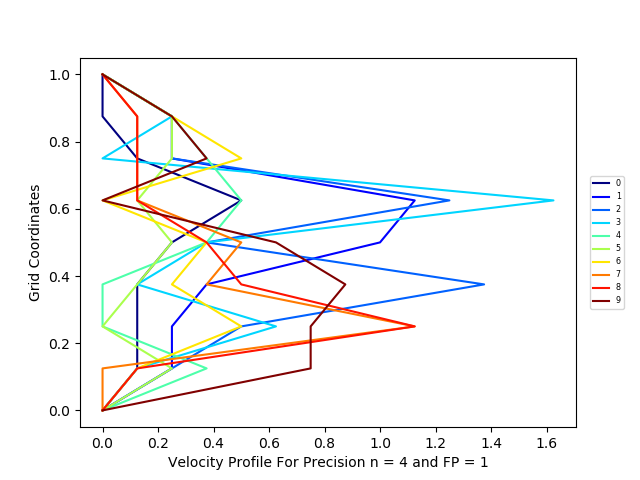}} \tabularnewline
\hline 
8 & {\footnotesize{}\includegraphics[scale=0.25]{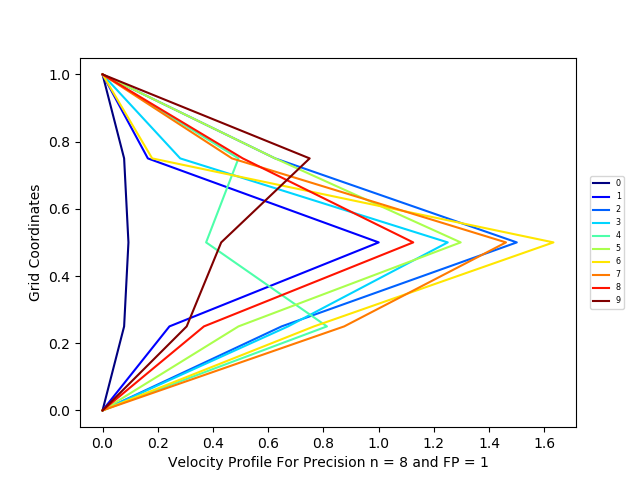}}
 & {\footnotesize{}\includegraphics[scale=0.25]{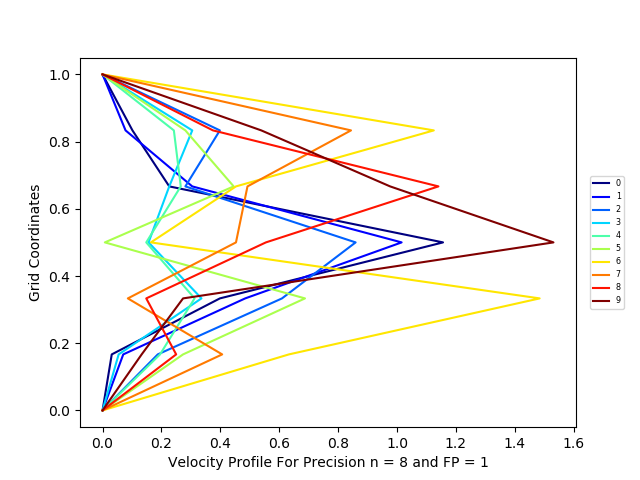}}
 & \tabularnewline
\hline 
double & {\footnotesize{}\includegraphics[scale=0.25]{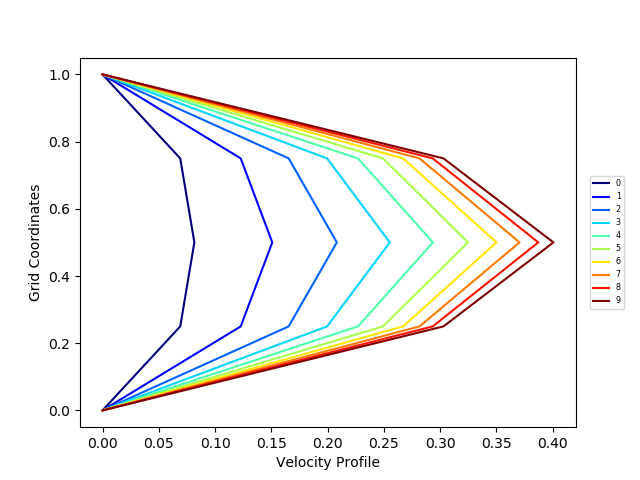}}
 & {\footnotesize{}\includegraphics[scale=0.25]{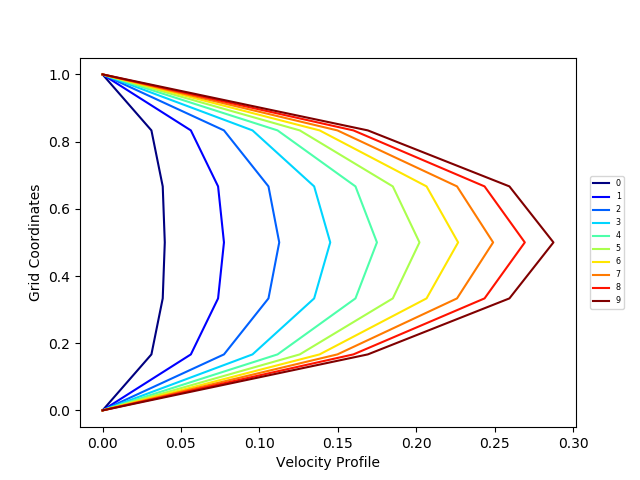}}
 & {\footnotesize{}\includegraphics[scale=0.25]{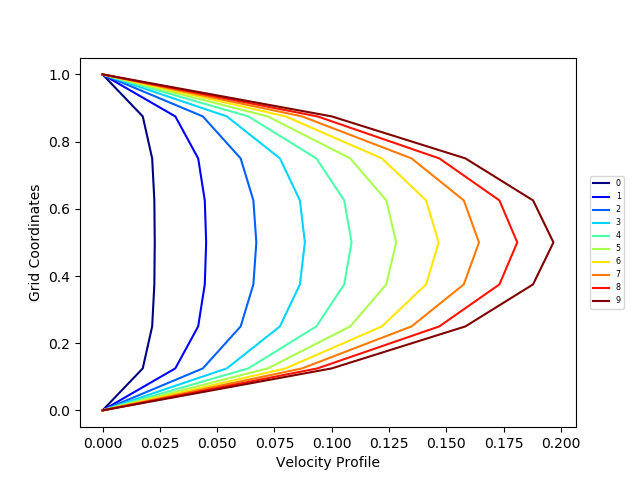}} \tabularnewline
% &  &  &  \tabularnewline
\hline  
\end{tabular}
\caption{{\large{}Solution profiles corresponding to the state with the lowest energy returned by the hardware for precision's 2, 4, 8 and double for mesh with sizes 5, 7 and 9.}}
\label{fig:lowest_energy}
\end{figure}

\begin{figure}[H]
\begin{tabular}{|>{\raggedright}m{0.15\paperwidth}|ccc|}
\hline 
$\#$ Grid Points & 5  & 7  & 9  \tabularnewline
\cline{1-1} 
$\#$ Precision &  &  &   \tabularnewline
\hline 
2 & {\footnotesize{}\includegraphics[scale=0.25]{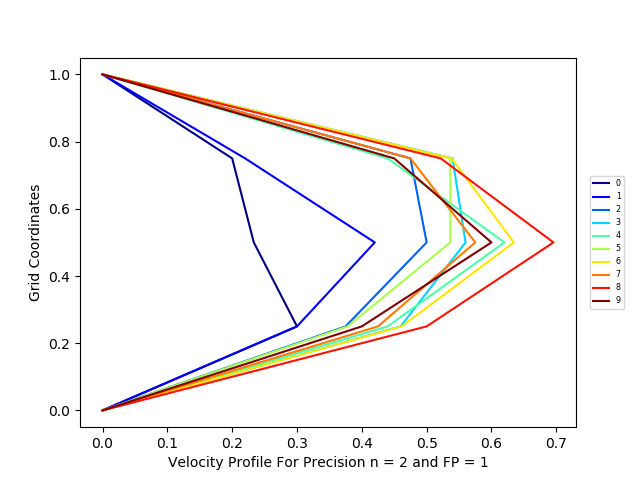}}
 & {\footnotesize{}\includegraphics[scale=0.25]{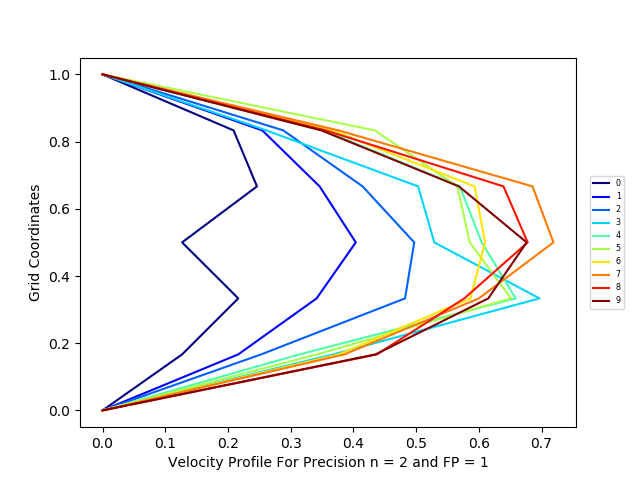}}
 & {\footnotesize{}\includegraphics[scale=0.25]{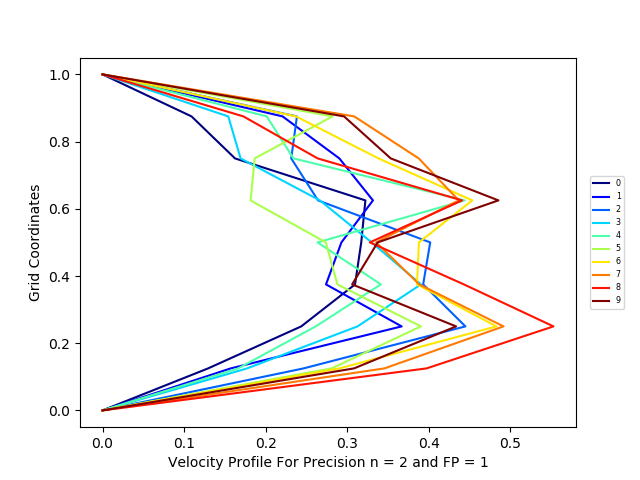}} \tabularnewline
\hline 
4 & {\footnotesize{}\includegraphics[scale=0.25]{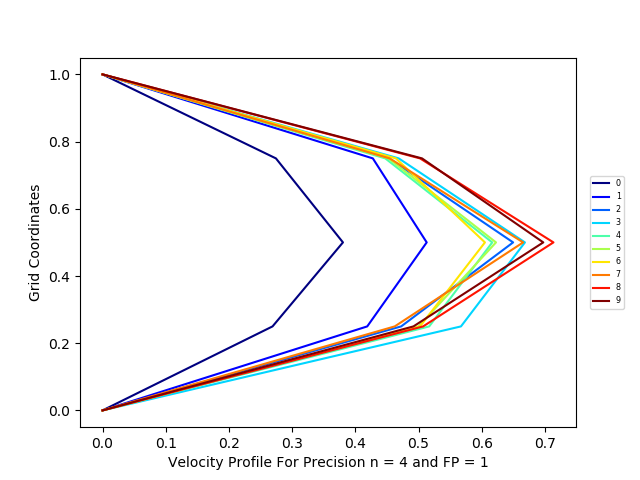}}
 & {\footnotesize{}\includegraphics[scale=0.25]{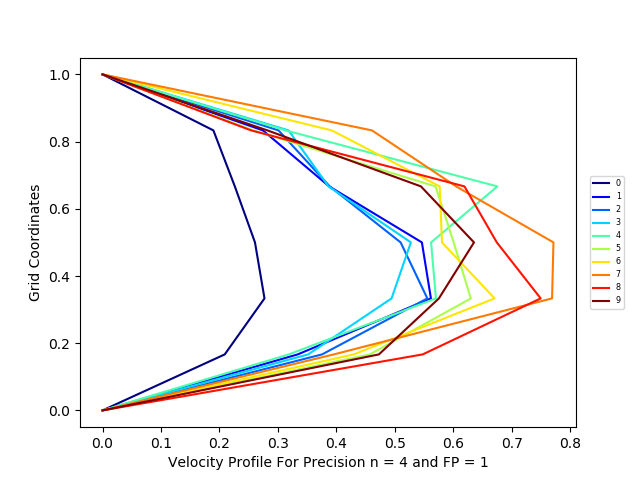}}
 & {\footnotesize{}\includegraphics[scale=0.25]{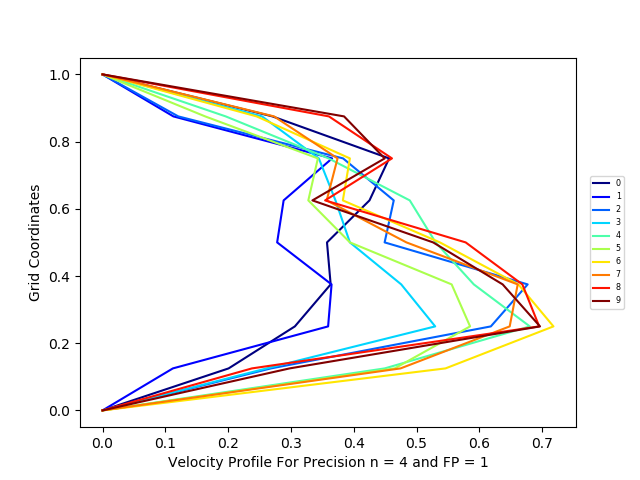}} \tabularnewline
\hline 
8 & {\footnotesize{}\includegraphics[scale=0.25]{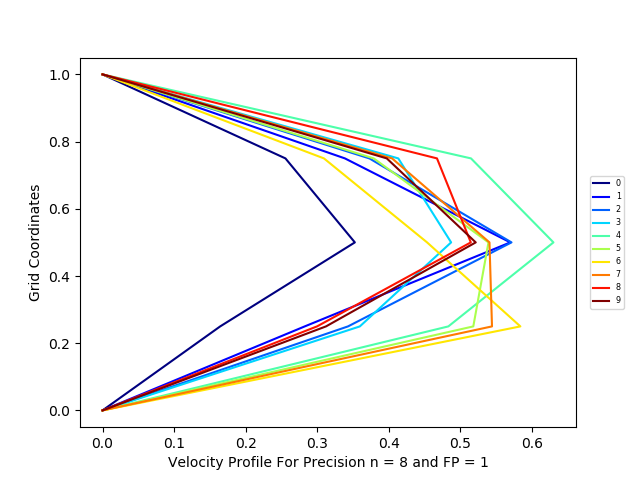}}
 & {\footnotesize{}\includegraphics[scale=0.25]{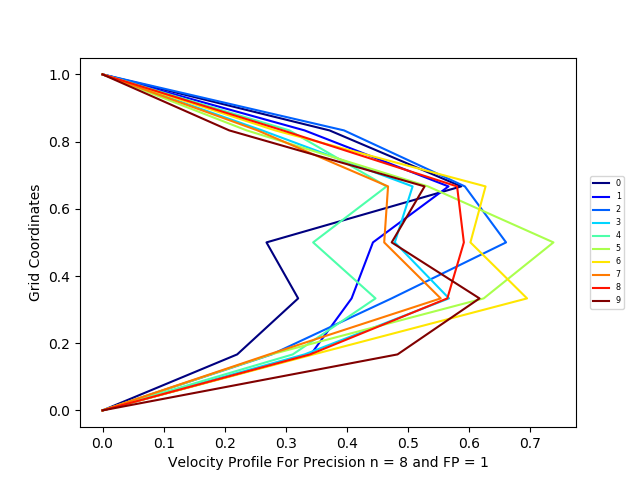}}
 & \tabularnewline
\hline 
double & {\footnotesize{}\includegraphics[scale=0.25]{./Figures/lowest_energy_soln_dw2x/Classical/classical_sol_nsteps_10_N_5}}
 & {\footnotesize{}\includegraphics[scale=0.25]{./Figures/lowest_energy_soln_dw2x/Classical/classical_sol_nsteps_10_N_7}}
 & {\footnotesize{}\includegraphics[scale=0.25]{./Figures/lowest_energy_soln_dw2x/Classical/classical_sol_nsteps_10_N_9}} \tabularnewline
\hline  
\end{tabular}
\caption{{\large{}Solution profiles corresponding to unweighted mean state for precision's 2, 4, 8 and double for mesh with sizes 5, 7 and 9.}}
\label{fig:umean_sol}
\end{figure}

\begin{figure}[H]
\begin{tabular}{|>{\raggedright}m{0.15\paperwidth}|ccc|}
\hline 
$\#$ Grid Points & 5  & 7  & 9  \tabularnewline
\cline{1-1} 
$\#$ Precision &  &  &   \tabularnewline
\hline 
2 & {\footnotesize{}\includegraphics[scale=0.25]{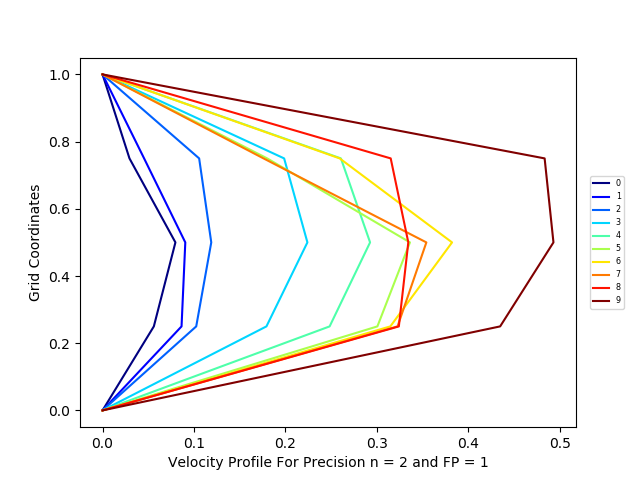}} 
 & {\footnotesize{}\includegraphics[scale=0.25]{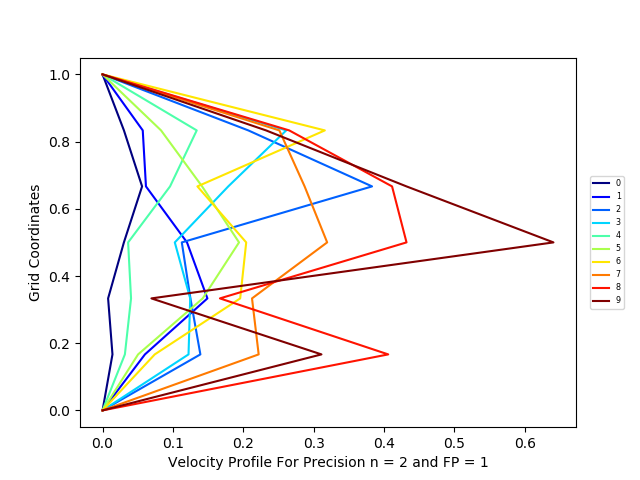}}
 & {\footnotesize{}\includegraphics[scale=0.25]{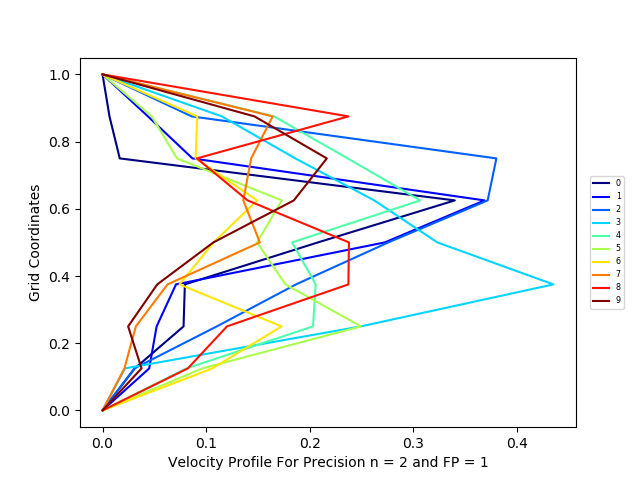}} \tabularnewline
\hline 
4 & {\footnotesize{}\includegraphics[scale=0.25]{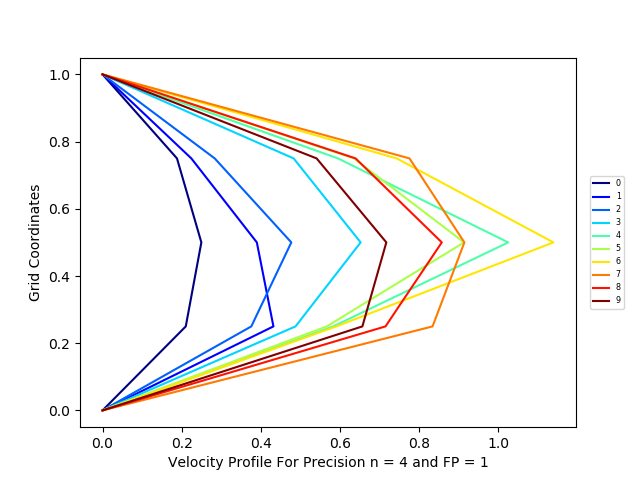}} 
 & {\footnotesize{}\includegraphics[scale=0.25]{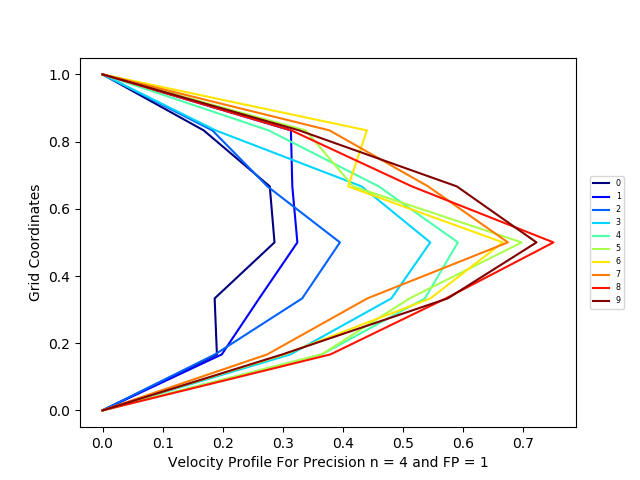}}
 & {\footnotesize{}\includegraphics[scale=0.25]{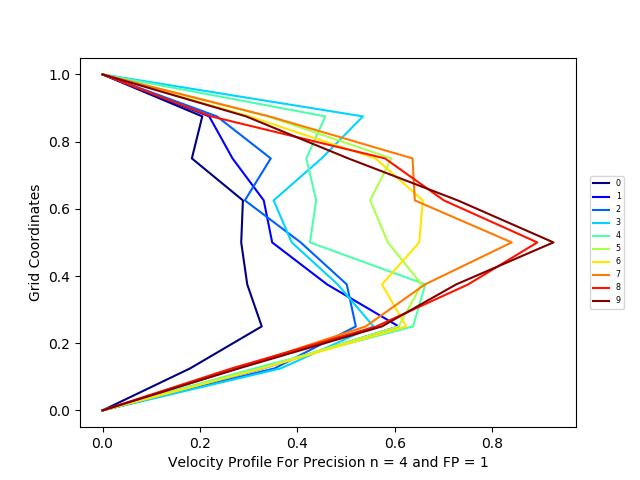}} \tabularnewline
\hline 
8 & {\footnotesize{}\includegraphics[scale=0.25]{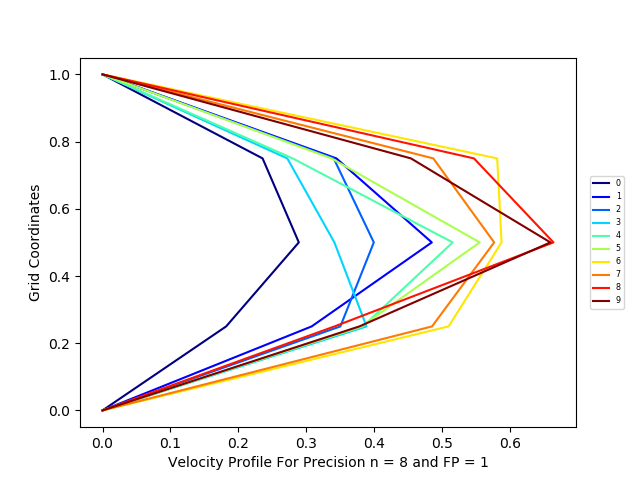}} 
 & {\footnotesize{}\includegraphics[scale=0.25]{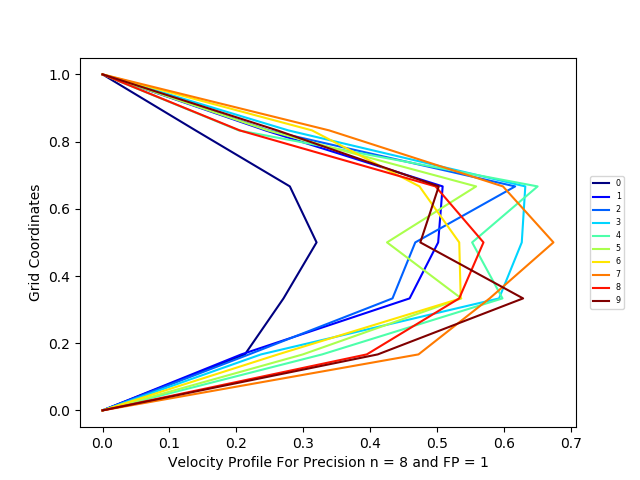}}
 & \tabularnewline
\hline 
double & {\footnotesize{}\includegraphics[scale=0.25]{./Figures/lowest_energy_soln_dw2x/Classical/classical_sol_nsteps_10_N_5}}
 & {\footnotesize{}\includegraphics[scale=0.25]{./Figures/lowest_energy_soln_dw2x/Classical/classical_sol_nsteps_10_N_7}}
 & {\footnotesize{}\includegraphics[scale=0.25]{./Figures/lowest_energy_soln_dw2x/Classical/classical_sol_nsteps_10_N_9}} \tabularnewline
\hline  
\end{tabular}
\caption{{\large{}Solution profiles corresponding to weighted mean state for precision's 2, 4, 8 and double for mesh with sizes 5, 7 and 9.}}
\label{fig:wmean_sol}
\end{figure}

Figure \ref{fig:lowest_energy}  shows the profiles of the lowest energy solutions obtained from DW2X\_3 solver for precision 2, 4 and 8. The results are shown for each of these three precisions for different grid resolutions - 5, 7 and 9 respectively. The last row shows the solution obtained using double-precision floating point arithmetic. For a mesh with 9 grid points, the logical problem solution size for precision 8 is 56 bits for which no embedding was found by the SAPI embedding utility. Figure \ref{fig:umean_sol} shows the results for the same combinations (precisions 2, 4, 8 and number of grid points 5, 7 and 9 as well as the classical solutions using double-precision floating point arithmetic)  but using an unweighted mean scheme, which means that the solution distributions obtained by SAPI are simply averaged, without any regard to their distributions (the number of times each solution vector is obtained and returned by SAPI). Figure \ref{fig:wmean_sol} shows the results for the same combination in figures \ref{fig:lowest_energy} and \ref{fig:umean_sol}, but now using a weighted mean of the solution vectors, weighted by the number of times they are returned by SAPI for each draw. The purpose of showing all three schemes in the figures \ref{fig:lowest_energy}, \ref{fig:umean_sol} and \ref{fig:wmean_sol} is to highlight the fact that no standardized techniques are available in the literature which prescribes how to extract meaningful solutions from the annealer for problems
of this nature as well as to study the sensitivity of the solutions to each of the three schemes described. For this purpose, all results are compared to the classical solution for the corresponding number of grid points. Also note that only 9 iterations of the solution process are shown for each solution to maintain clarity in presentation and also due to limitations of resources in terms of allocation time on the machine and resources available for data analysis and postprocessing. However, this exercise is still useful to highlight the broad sensitivities of the solution space. We highlight several observations from figures \ref{fig:lowest_energy} below.

\begin{itemize}
\item The lowest energy solution state returned by the solver for lower precision's are highly inaccurate.  For precision 2, the solution yields a straight line, which is unrealistic.  
\item Unlike the classical solutions, which shows a systematic monotonic progression over each time step, the least energy solutions show a much more chaotic oscillatory behavior with no systematic progression over iterations steps. 
\item Only the solution with precision 8 resembles the shape of the classical profile qualitatively, although very crudely. 
\item Quantitatively, the peak values of the profile (expected to be parabolic from the classical solution) are higher by an order of magnitude. More details on the quantitative comparisons are provided in section \ref{sec:error_analysis}.
\end{itemize}

While the lowest energy solutions seem to be  provide very crude approximations, the schemes with weighted and unweighted means perform better in comparison. We summarize the key observations from figure \ref{fig:umean_sol} next. 

\begin{itemize}
\item Most of the solutions show a systematic progress of the solution with the successive iterations. 
\item For all precision values, lower number of grid points (5 and 7) qualitatively capture the parabolic solution profile, although solution profiles start to degrade with the increase of precision and number of grid points.
\item With 9 grid points, the shape of the profile appears bimodal and is unrealistic.
\item The peak value of the profile still deviates from the classical solution, however, the overall solution quality is better than the least energy formulation, in terms of error magnitude.
\end{itemize} 

Hence it appears that using all the states than a single state based on the least energy can improve the quality of solution, but it is also dependent on the precision and grid size. To understand the sensitivity of the solutions when all states are used and weighted by their number of occurrences to obtain the mean is studied next. The following key observations are derived from figure \ref{fig:wmean_sol}. 
\begin{itemize}
\item Solutions do progress more systematically like the unweighted scheme, at least for 5 and 7 grid points.
\item For 5 grid points, solutions are quite close to the classical solutions across the iterations, and the solution quality improves with precision.
\item For higher number of grid points, the solutions seem to improve with higher precision, which is an opposite behavior compared to the unweighted cases.
\item The errors appear to be reduced with increase in precision and number of grid points.
\end{itemize}

\subsection{Solution Distribution at Domain Center}
\begin{figure}[H]
\begin{tabular}{|>{\raggedright}m{0.06\paperwidth}|cccc|}
\hline 
Precision & 2 & 4 & 6 & 8\tabularnewline
\cline{1-1} 
$\#$ Grid Points &  &  &  & \tabularnewline
\hline 
\multirow{3}{0.06\paperwidth}{5} &  &  &  &  \tabularnewline
 & {\footnotesize{}\includegraphics[scale=0.22]{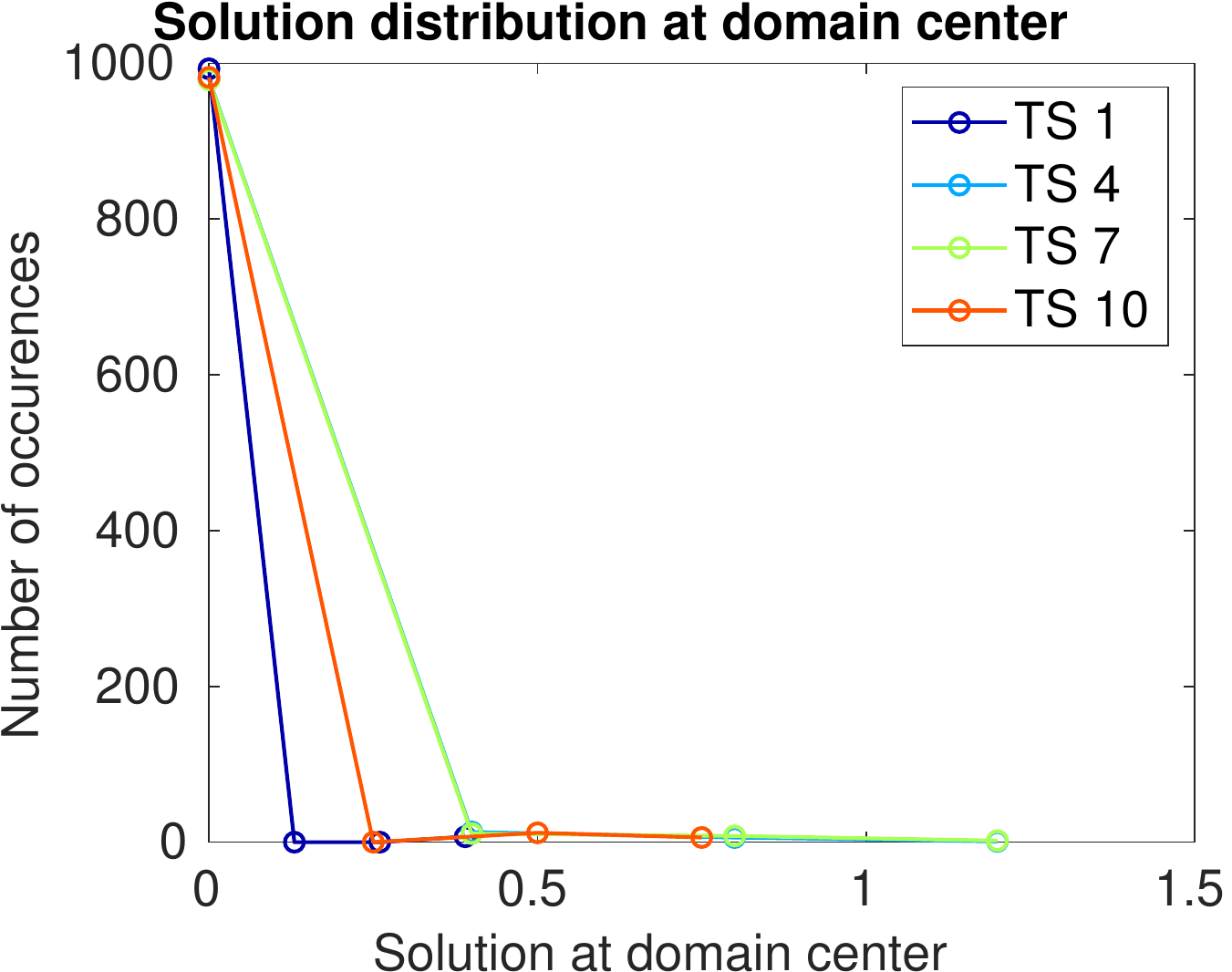}} & {\footnotesize{}\includegraphics[scale=0.22]{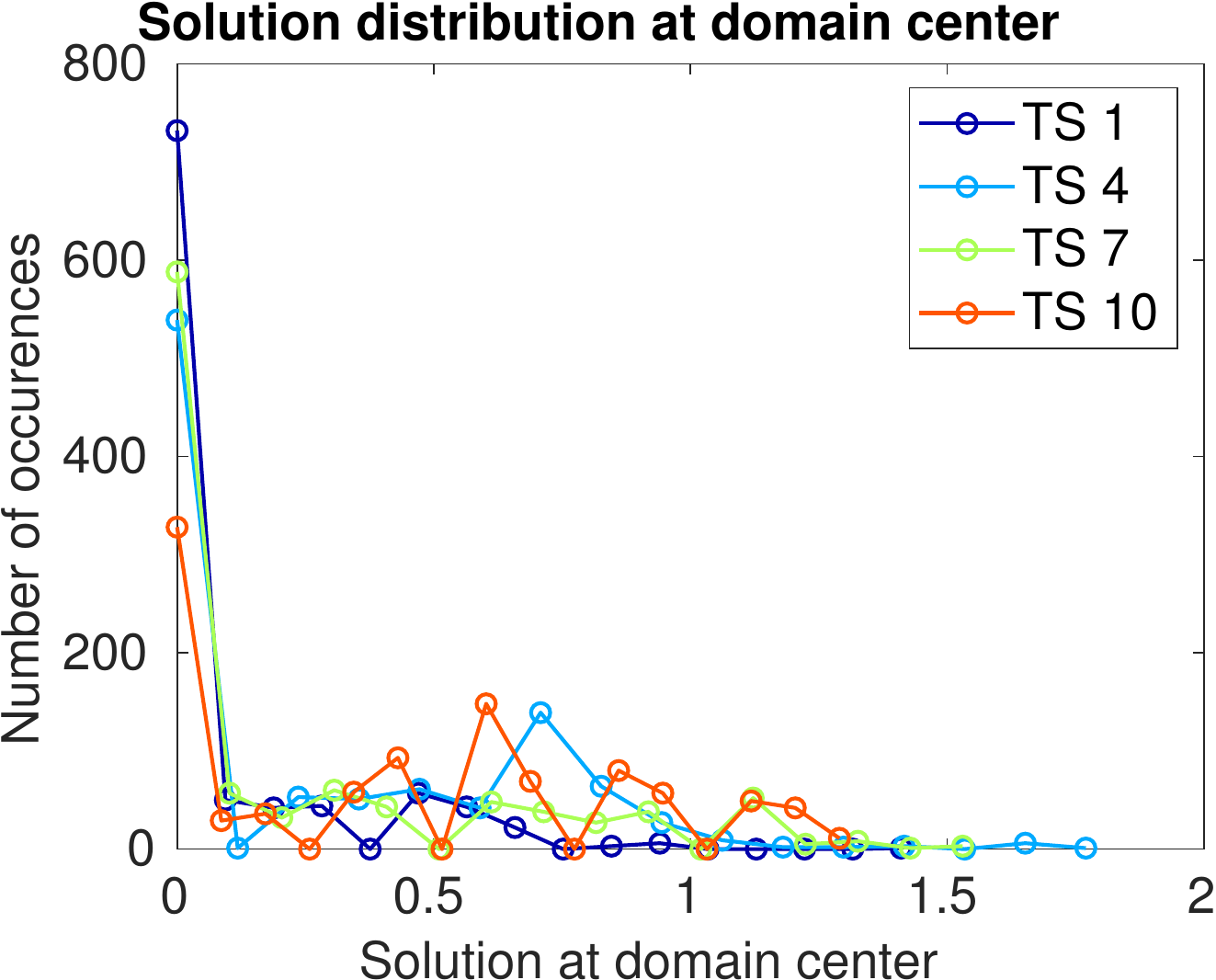}} & {\footnotesize{}\includegraphics[scale=0.22]{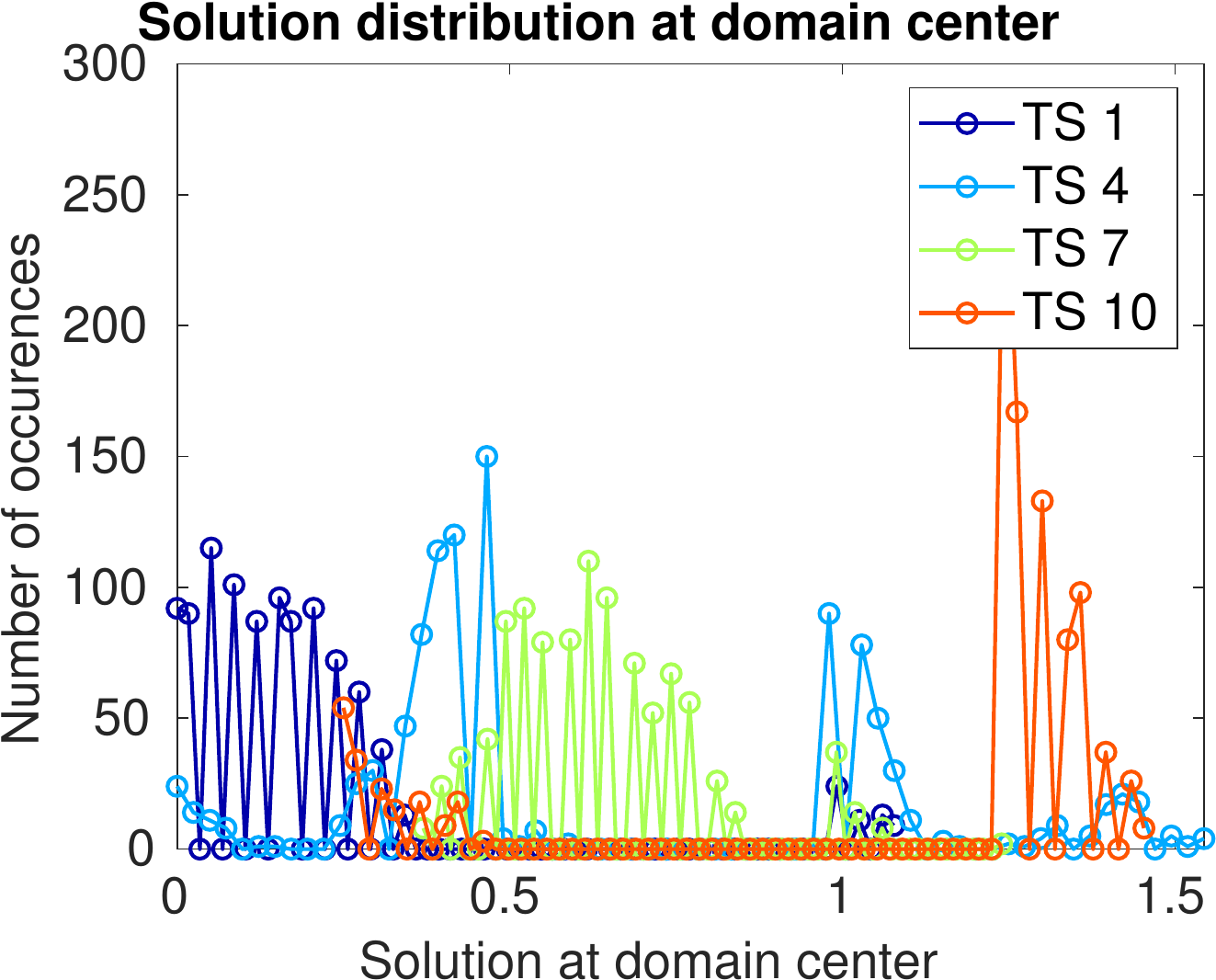}} & {\footnotesize{}\includegraphics[scale=0.22]{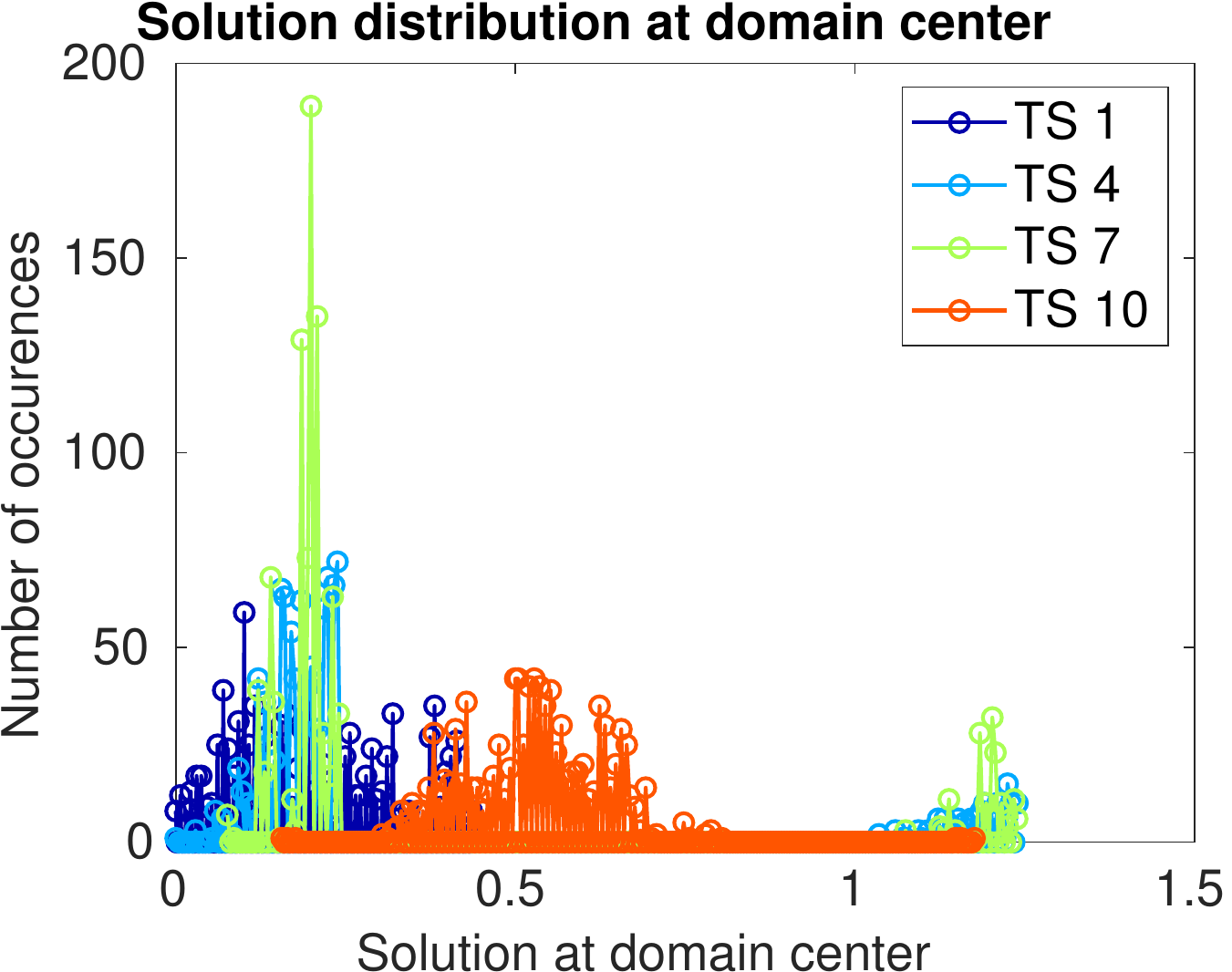}}\tabularnewline
 &  &  &  & \tabularnewline
\hline 
\multirow{3}{0.06\paperwidth}{7} &  &  &  & \tabularnewline
 & {\footnotesize{}\includegraphics[scale=0.22]{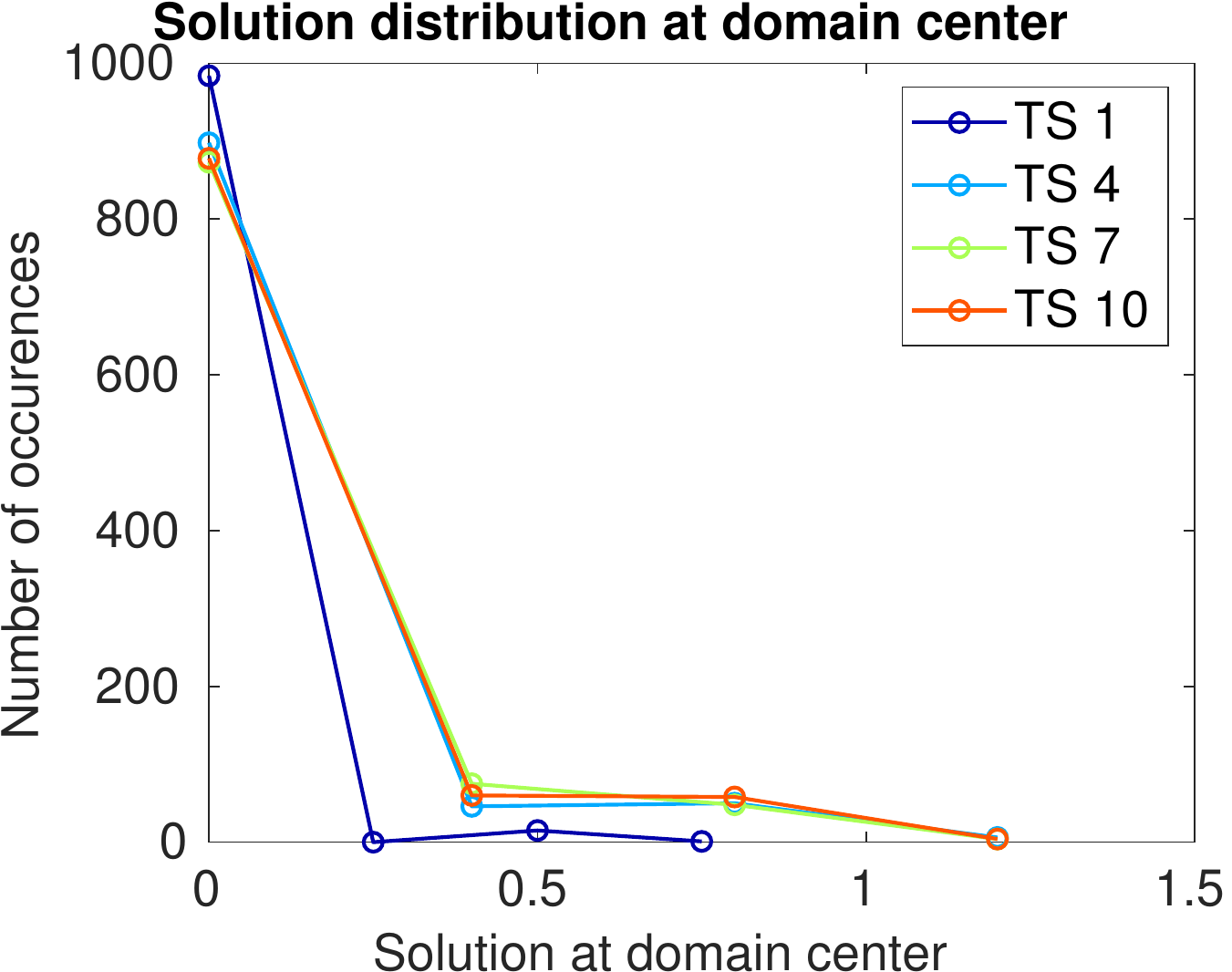}} & {\footnotesize{}\includegraphics[scale=0.22]{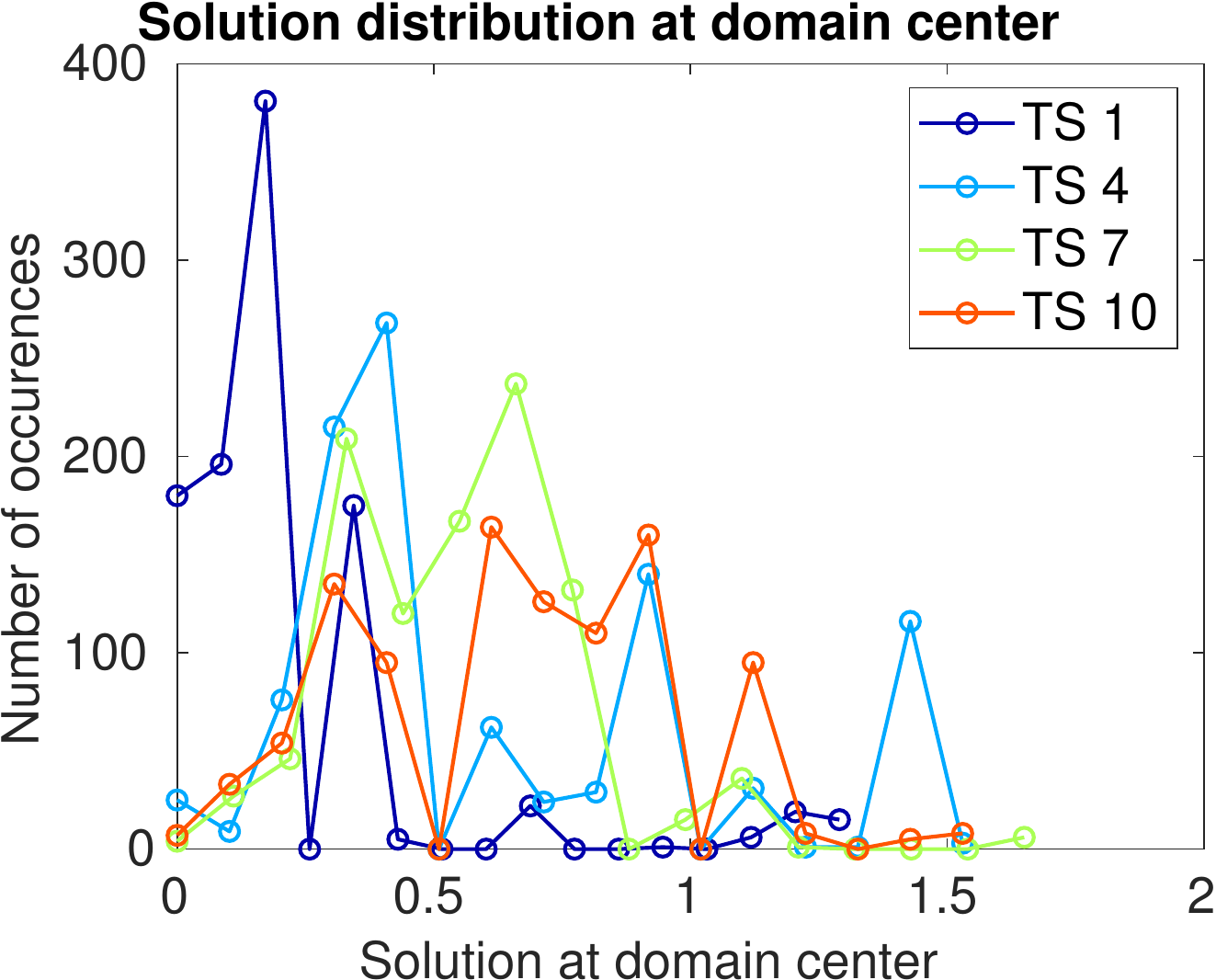}} & {\footnotesize{}\includegraphics[scale=0.22]{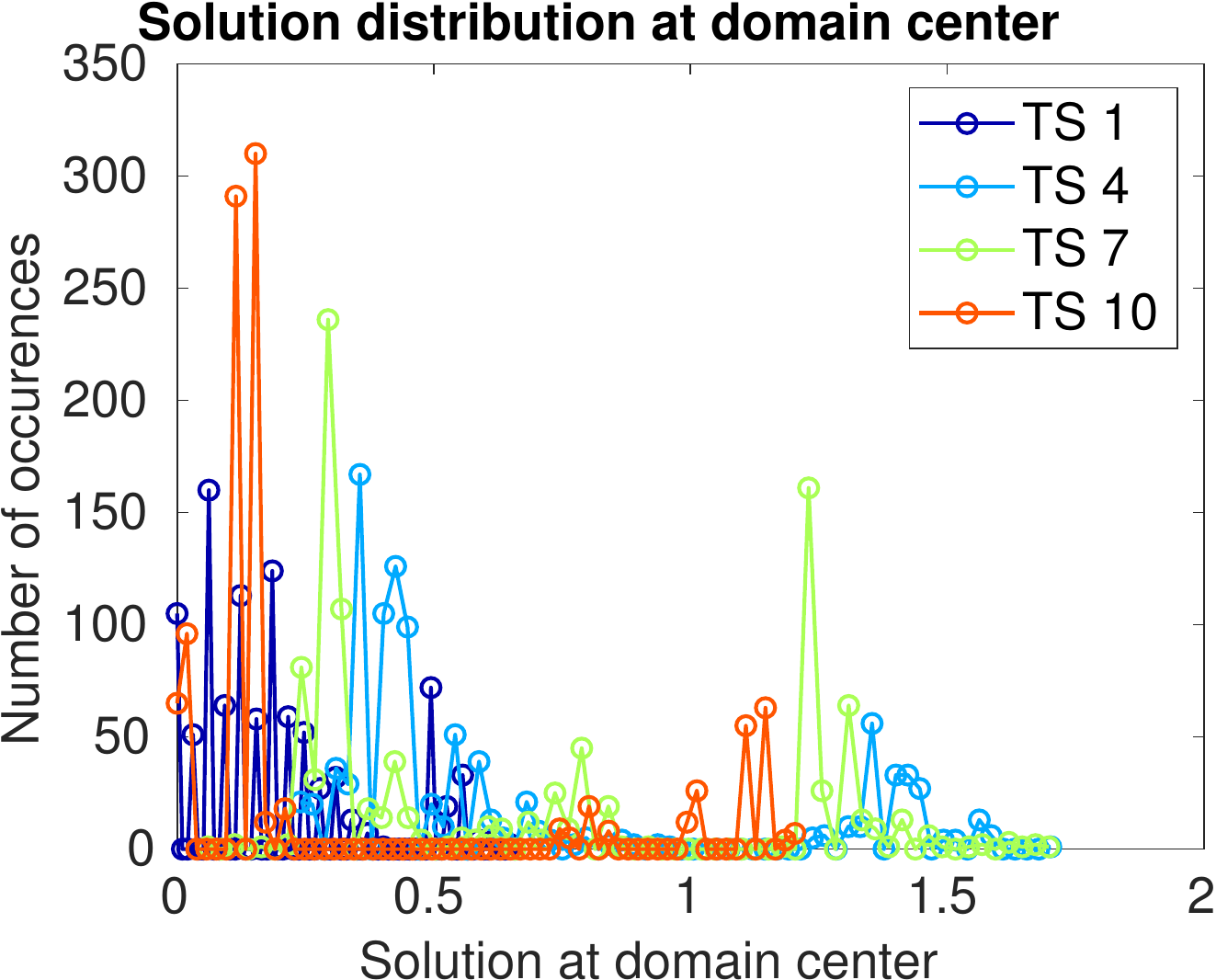}} & {\footnotesize{}\includegraphics[scale=0.22]{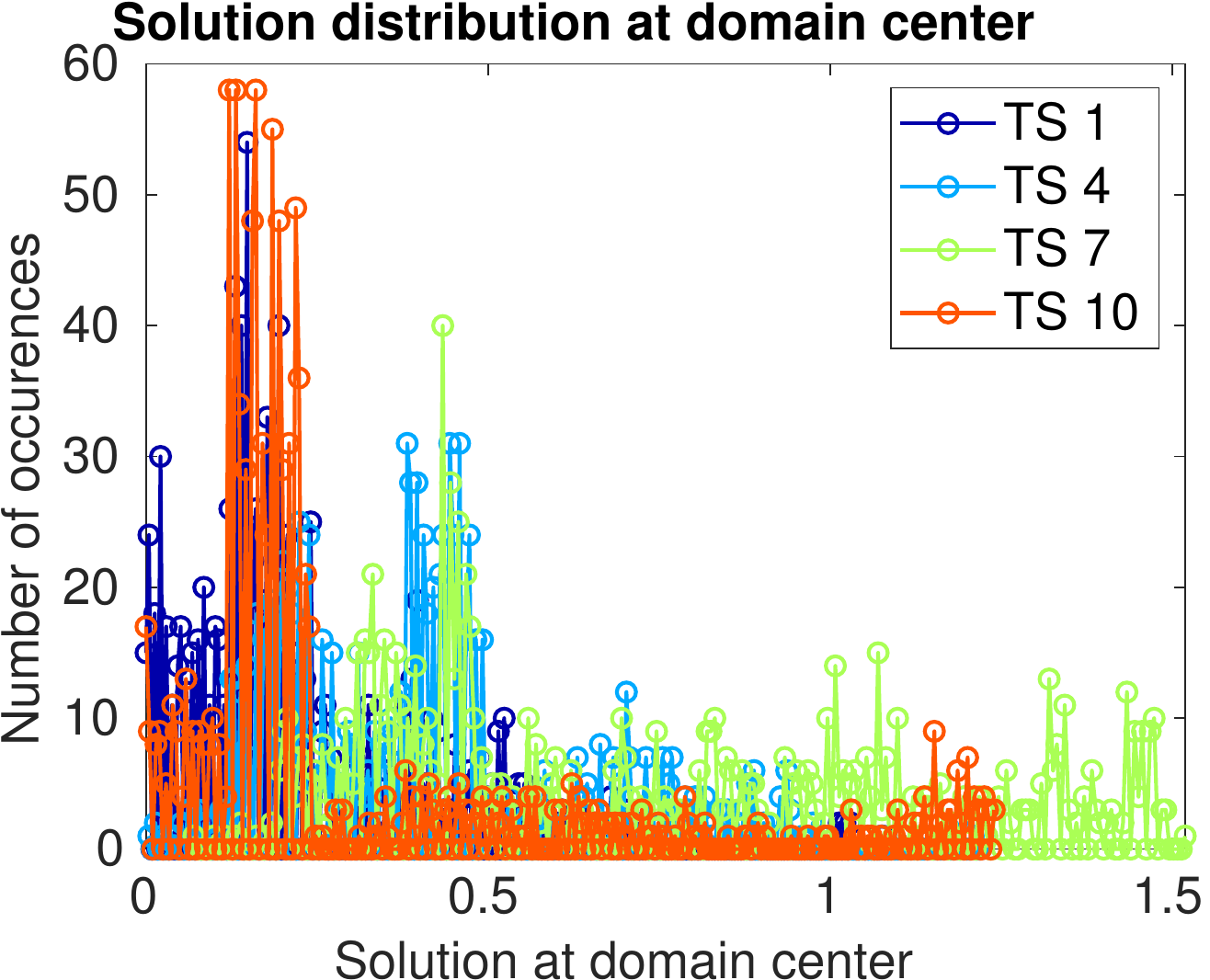}}\tabularnewline
 &  &  &  & \tabularnewline
\hline 
\multirow{3}{0.06\paperwidth}{9} &  &  &  & \tabularnewline
 & {\footnotesize{}\includegraphics[scale=0.18]{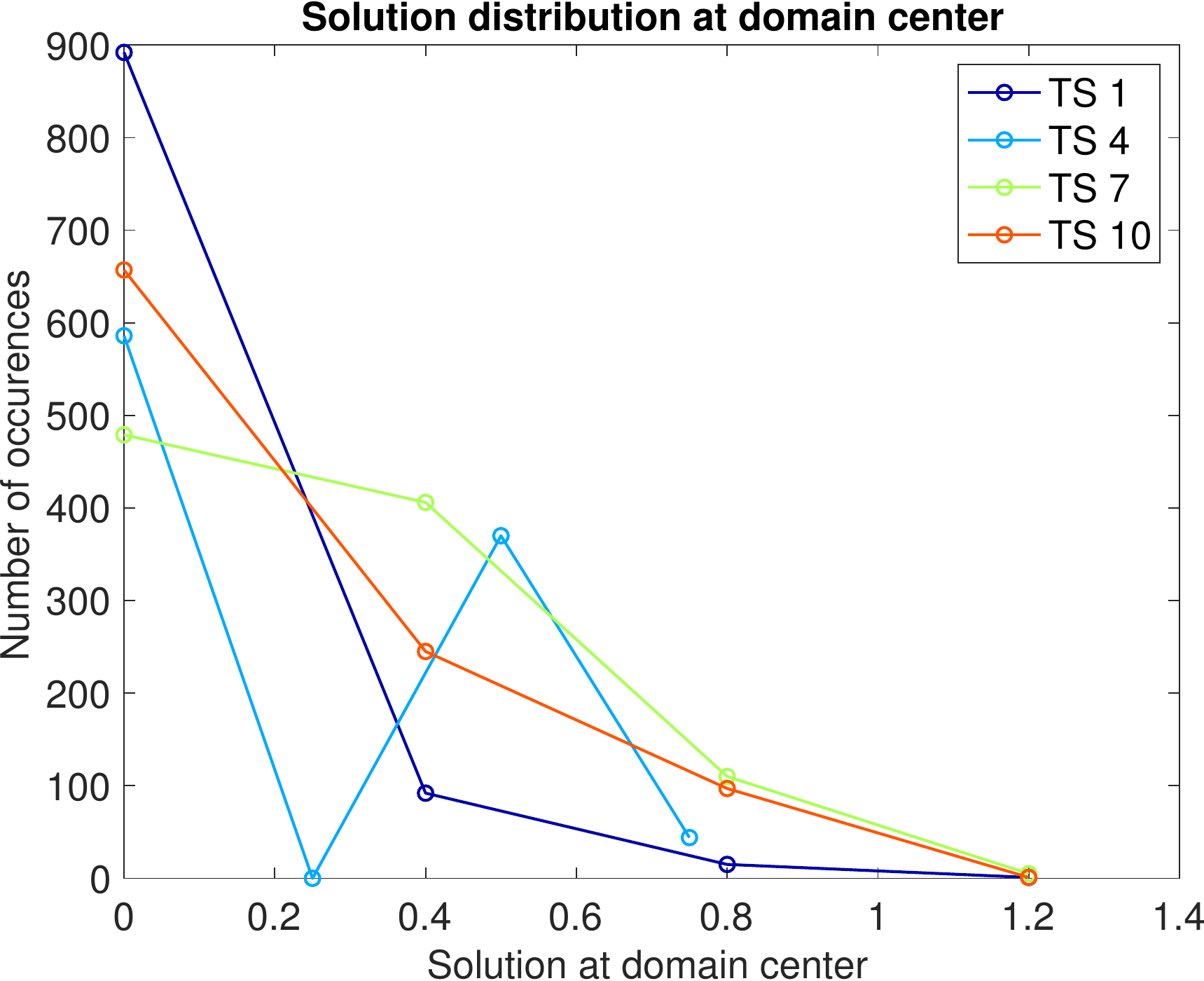}} & {\footnotesize{}\includegraphics[scale=0.18]{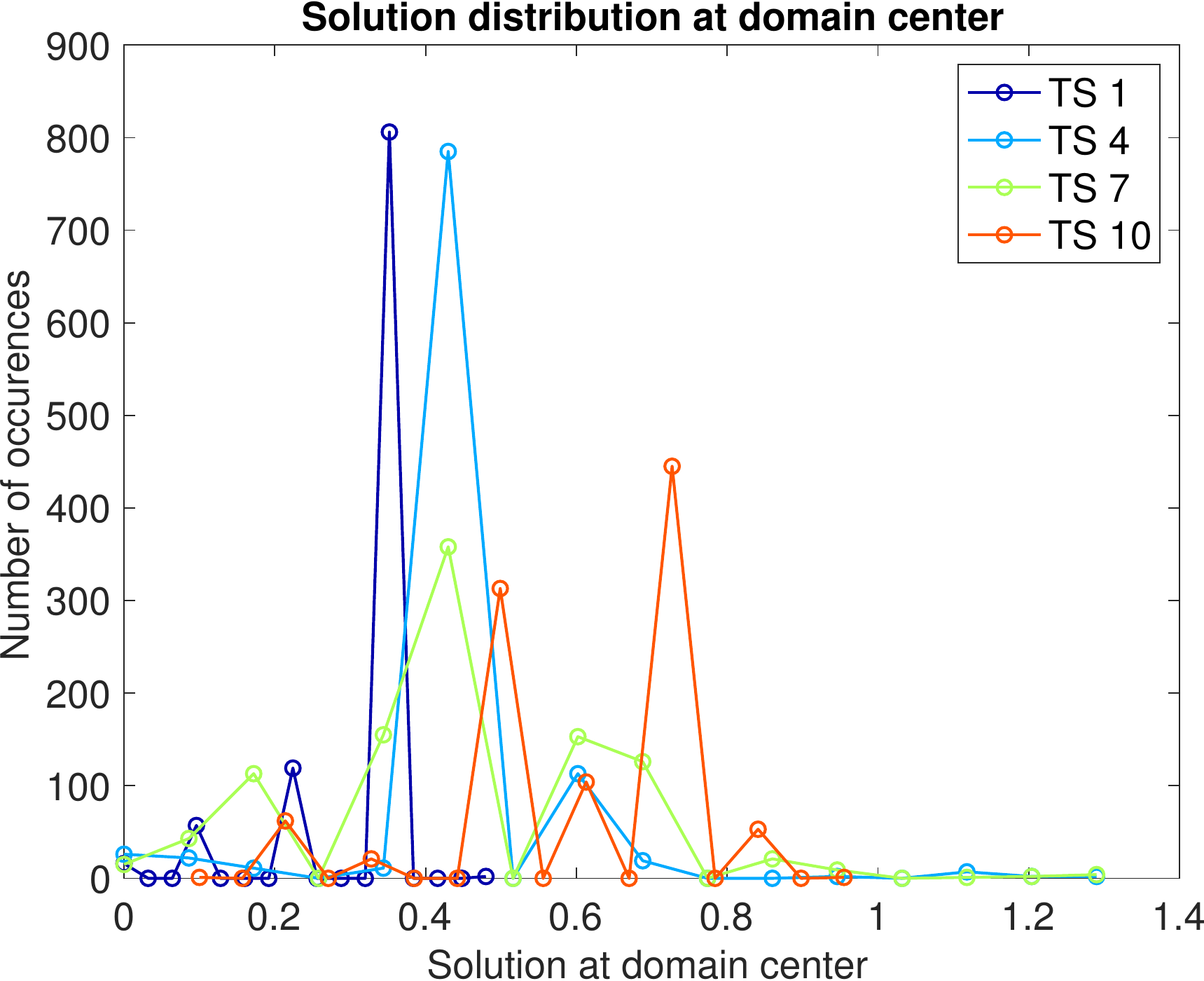}} &   {\footnotesize{}
 \includegraphics[scale=0.18]{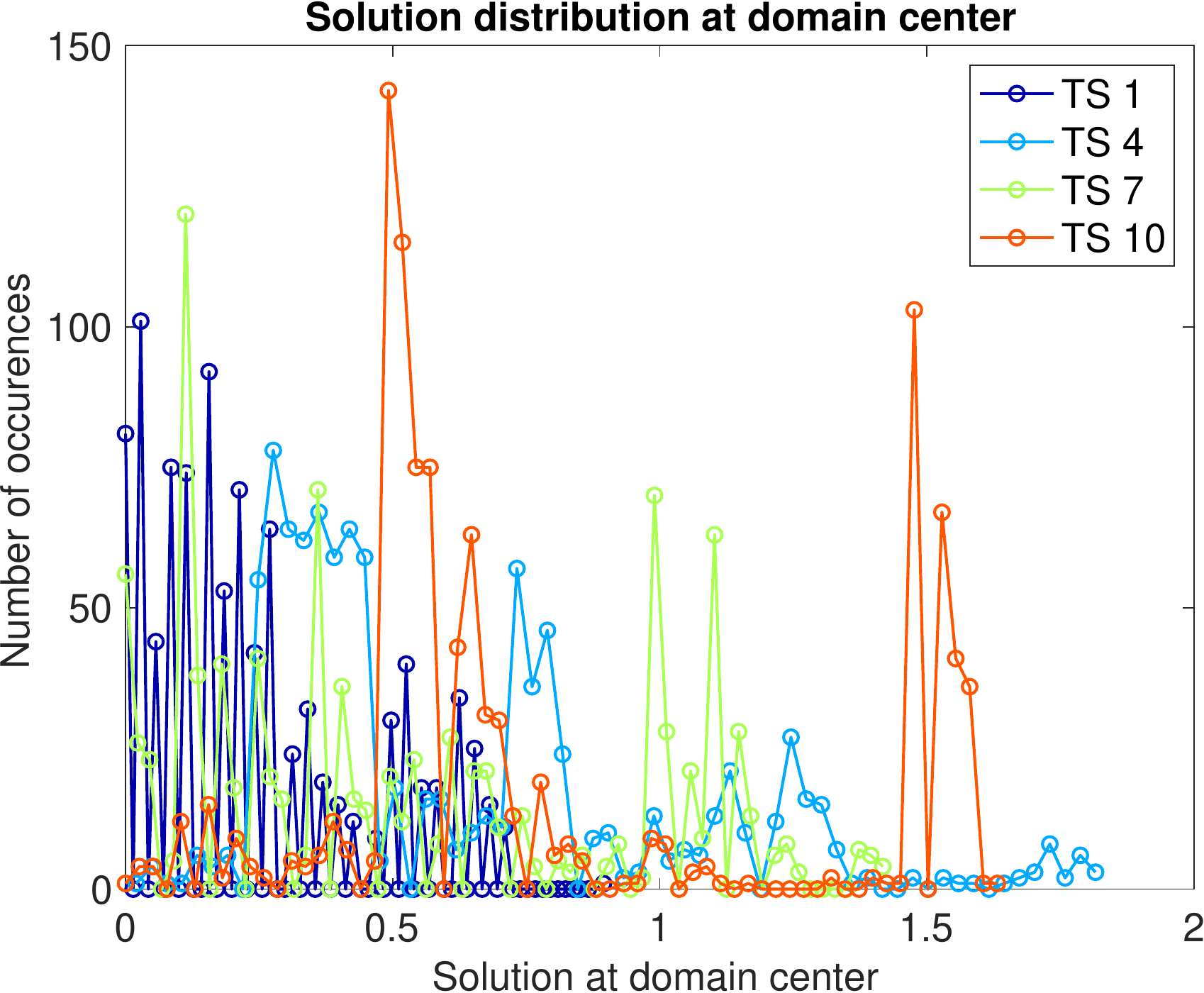}} & NA \tabularnewline
 &  &  &  & \tabularnewline
\hline 
\end{tabular}
\caption{{\large{}Solution Distribution at Domain Center}}
\label{fig:sol_distribution}
\end{figure}

To explain the difference in behavior between the schemes using unweighted and weighted means, the actual distributions of solutions at the domain center are plotted in figure \ref{fig:sol_distribution}. It shows the aforementioned distributions (actual numbers of occurrences) for precision values of 2, 4, 6 and 8, each for 5, 7 and 9 grid points. For clarity of presentation, only 4 different iteration steps are plotted (1,4, 7 and 10). The key observations from figure \ref{fig:sol_distribution} are: 

\begin{itemize}
\item For lower precision values, the number of occurrences of the state which corresponds to the real value zero occurs  the highest number of times. With increase in the number of iterations and the number of grid points, non zero solutions start to appear a finite number of times. Thus, taking a simple unweighted mean will counter the effects of high frequencies of zero solutions. However, taking a weighted mean will push the final solution towards zero. This explains why taking an unweighted mean yields a better solution at lower precision values.

\item At higher precision values, the non zero solutions start having more frequencies as expected. As the precision is increased further, the distributions also progressed finitely in time. As noticed in the figure, the time iterations of the distributions can be identified as separate clusters.

\item At higher precision, taking a weighted mean seems to be a better strategy as that approach will represent the finite frequencies of non zero solutions across the solution space.
\end{itemize}

\subsection{Error Analysis}
\label{sec:error_analysis}
\begin{figure}
     \begin{center}
         \includegraphics[width=\linewidth]%
             {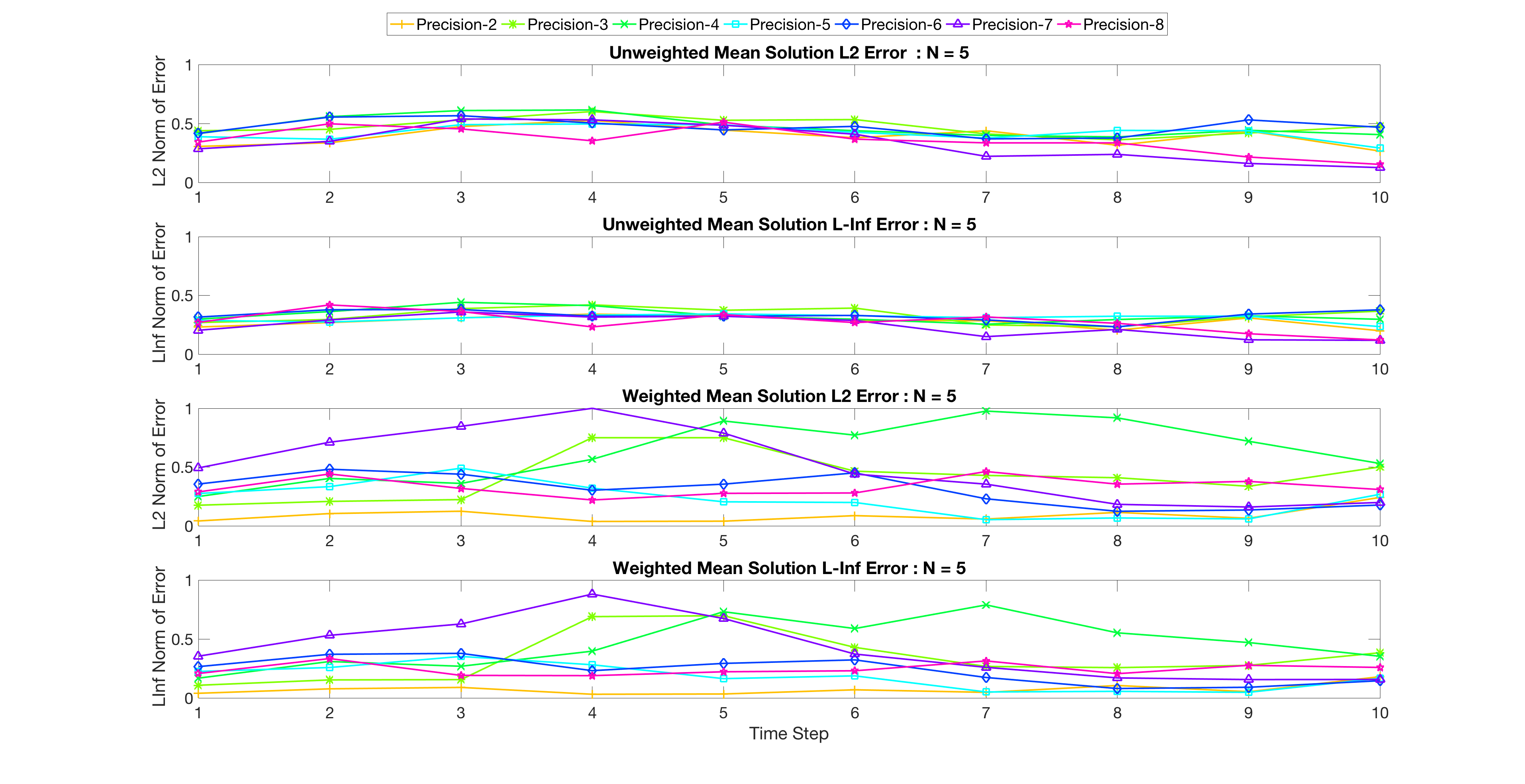}
     \end{center}
     \caption{{\large{} $L_{2}$ and ${L_{\infty}}$ errors of the unweighted and weighted mean solutions for a problem with 5 grid points. The precision varies from 2 to 8 and the error is computed for each time step. }}
     \label{fig:error_n5}
 \end{figure}

\begin{figure}
     \begin{center}
         \includegraphics[width=\linewidth]%
             {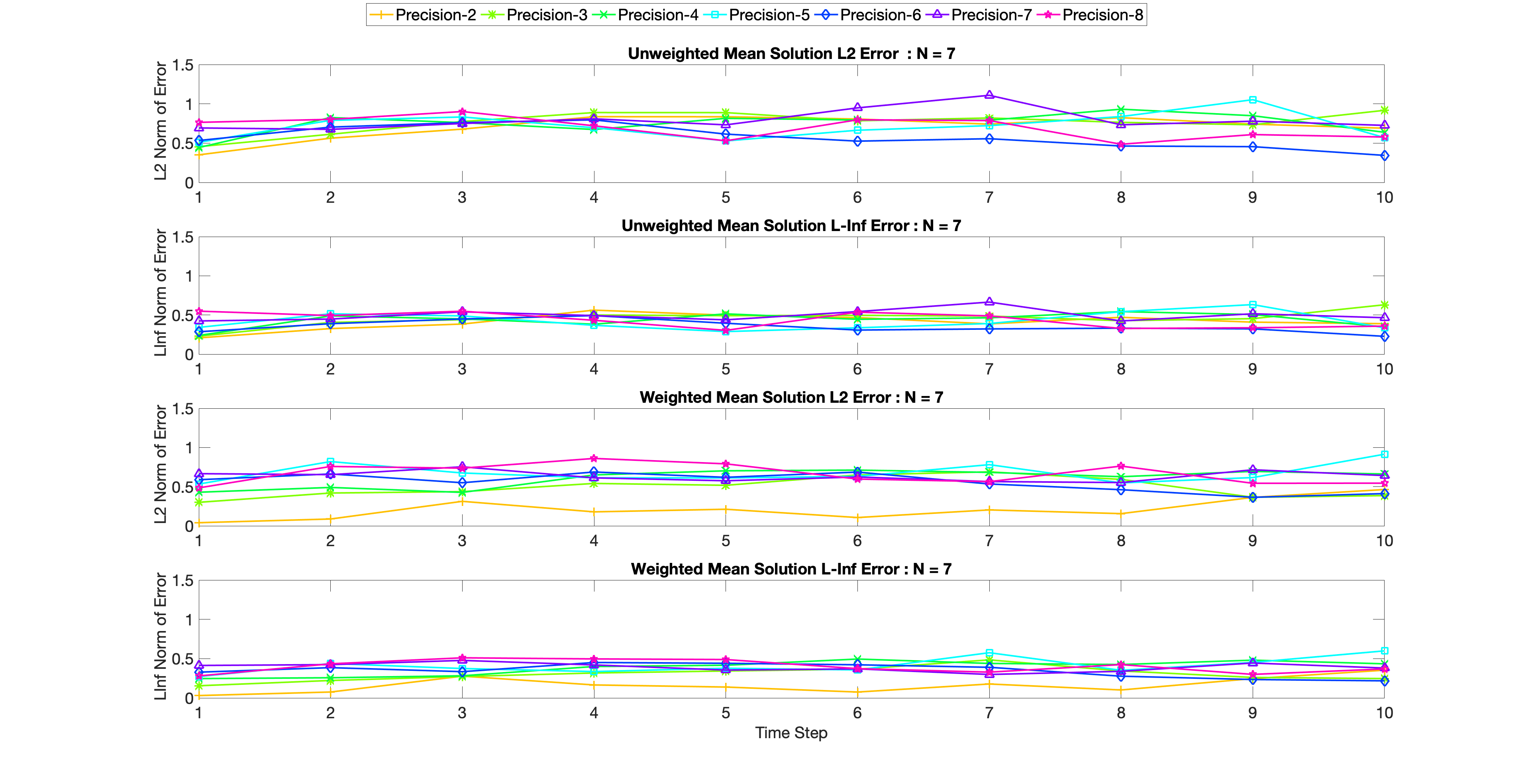}
     \end{center}
      \caption{{\large{} $L_{2}$ and ${L_{\infty}}$ errors of the unweighted and weighted mean solutions for a problem with 7 grid points. The precision varies from 2 to 8 and the error is computed for each time step. }}
      \label{fig:error_n7}
 \end{figure}

\begin{figure}
     \begin{center}
         \includegraphics[width=\linewidth]%
             {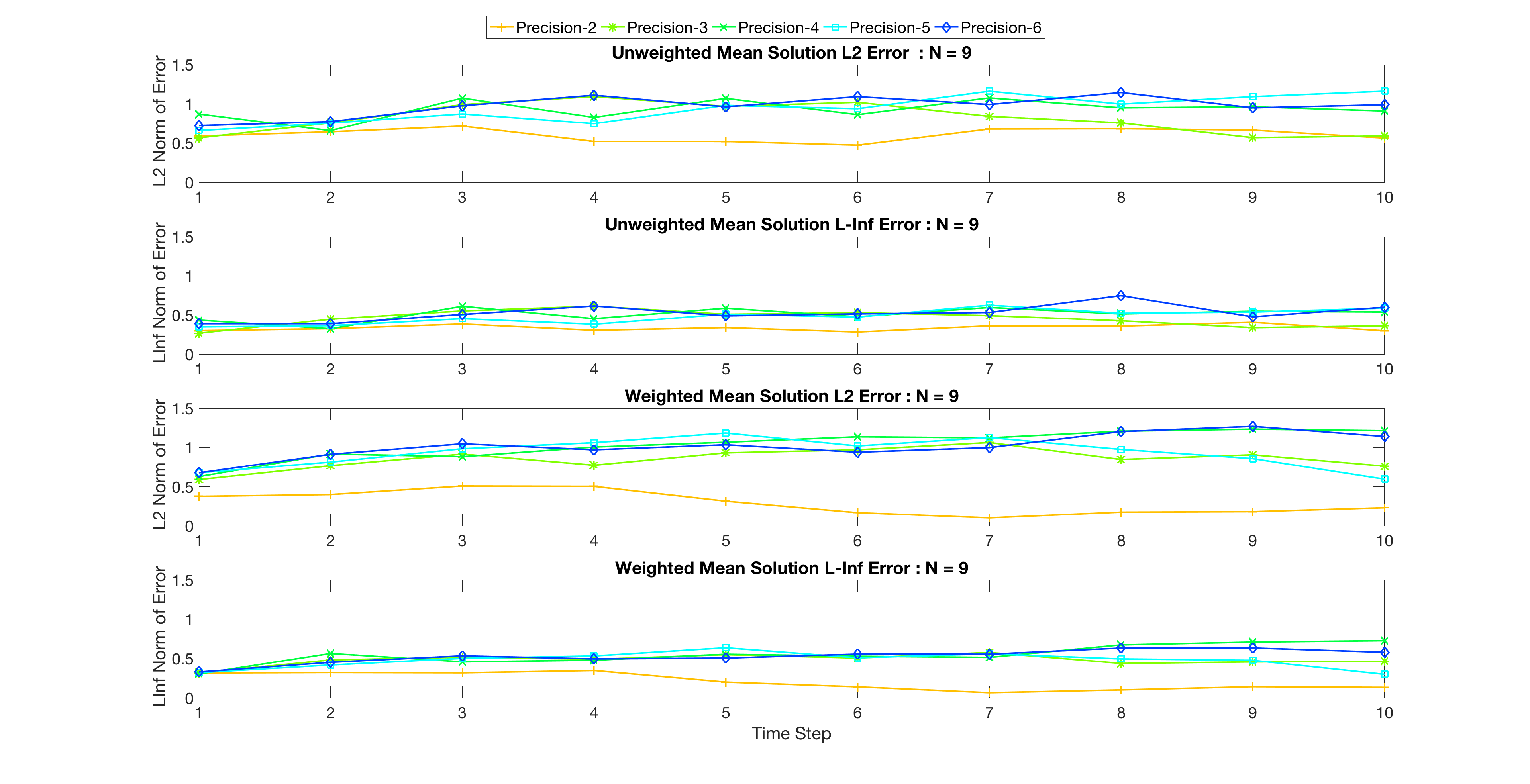}
     \end{center}
      \caption{{\large{} $L_{2}$ and ${L_{\infty}}$ errors of the unweighted and weighted mean solutions for a problem with 9 grid points. The precision varies from 2 to 8 and the error is computed for each time step. }}
      \label{fig:error_n9}
 \end{figure}

Figures \ref{fig:error_n5}, \ref{fig:error_n7} and \ref{fig:error_n9} summarizes the solutions by plotting the $L_2$ and $L_{\infty}$ errors for 5, 7 and 9 grid points, respectively. Each of these figures show both types of errors for the unweighted and weighted schemes. The $L_2$ error computes the Euclidean norm of the difference between the quantum solution and that of the classical solution. The $L_{\infty}$ error effectively computes the difference between the maximum value of the solution vector for the classical and the quantum solutions. Hence  the $L_2$ error represents the error for the entire solution space, and the $L_{\infty}$ error represents the error for the solution maximum, which represents the highest velocity within the domain, i.e., the center point of the domain. The errors are shown for all precision values ( 2 to 8) in figures \ref{fig:error_n5}, \ref{fig:error_n7} and \ref{fig:error_n9}.  We summarize the key observations from these figures below:
\begin{itemize}
\item For 5 grid points, the difference among different precision values are much lower for the unweighted scheme compared to the weighted scheme for both $L_2$ and $L_{\infty}$ errors.

\item For 5 grid points, both the $L_2$ and $L_{\infty}$ errors are lower for higher precision values for unweighted solutions. Moreover, the errors seem to decrease with iterations, indicating potential convergence. The weighted mean solutions do not exhibit any clear pattern like this.

\item For 7 and 9 grid points, the divergence among the weighted solutions seem to be larger compared to the unweighted scheme. However, errors seem to increase with higher precision and the lowest precision value seems to show the lowest  $L_2$ and $L_{\infty}$ errors. No clear pattern of convergence is found either.
\end{itemize}

%==============
%Lowest energy, #grid points {5,7,9}, precision {8,4,2}, fp 1  

%==============
%Unweighted mean energy, #grid points {5,7,9}, precision {8,4,2}, fp 1  

%==============
%Weighted mean energy, #grid points {5,7,9}, precision {8,4,2}, fp 1  

%============
% Error Plots
 
\section{Conclusions} \label{conclusions}
In this work, we have demonstrated the potential, shortcomings and sensitivities of solving linear systems obtained from a reduced version of the Navier Stokes (NS) equation on the DWave quantum computer. We have chosen a very simple application, namely the one dimensional channel flow problem, where the classical solution is a well established parabolic flow. Using this reference solution we have explored how the system of linear equations obtained from the discretized form of the NS equations can be converted to a form amenable to the DWave Quantum annealer. This conversion involves a  fixed point arithmetic based conversion of decimal to binary variables and posing the system in a least square formulation. This allows the construction of a quadratic unconstrained binary optimization (QUBO) problem, that is solvable by the DWave utilities. However, in this work, we have only used the default embedding of the DWave annealer (called a Chimera graph) which is organized as a lattice of unit blocks of qubits, where each block has eight qubits configured as a four-node bipartite graph. A technique called chaining allows one to use more optimized embeddings specific to the problem, and this is beyond the scope of the current work. In this work we have instead focused on the following questions:
\begin{itemize}

\item How to get extract the `correct' solution from the distribution of binary solution vectors?

\item How do the different strategies for obtaining this solution compare to each other?

\item What are the sensitivities of the obtained solutions to the grid resolution and precision value used in the fixed point arithmetic?
\end{itemize}

To answer these question we have plotted solution distributions at the domain center, compared three different methods (least energy, simple unweighted means and weighed means, where the weights are number of occurrences of each solution) to extract the quantum solutions against the classical solution, and performed error analysis to demonstrate the sensitivities with precision values and grid resolutions. We can answer the aforementioned questions as follows:
\begin{itemize}

\item There is no unique method to extract the correct solution.

\item The least energy solution is the poorest. The weighted and unweighted means that involve all obtained solutions in each iterations provide better solutions consistently. 

\item At lower precision values, the unweighted means perform better. At higher precision values weighting the solutions with their numbers of occurrences is a better approach. Interestingly, increase of grid resolution and higher precision does not automatically allow more accurate solutions because it increases the noise of the system significantly. 
\end{itemize}

To the best of our knowledge of the literature, this work is one of the first which demonstrates a simple workflow to reduce the Navier Stokes equations in a form readily solvable by the DWave machine (or quantum computers of the annealer type). This workflow can be followed to solve several problems spanning multiple disciplines ranging from solid to fluid mechanics and engineering, essentially any problem that can be reduced to linear systems of equations. The sensitivity of the workflow to the multiple parameters of the QUBO type solution process is perhaps of more interest than the accuracy of the solution instead. In future works, we will further explore the effects of more optimized embedding schemes. Moreover, the inherent nature of the DWave machine, which provides a large number of distributions at every step of the solution instead of one unique solution might allow one to explore Monte Carlo type approaches to solving linear systems of practical interest.

%%%%%%%%%%%%%%%%%%%%%%%%%%%%%%%%%%%%%%%%%%%%%%%%%%%%%%%%%%%%%
\section{Acknowledgments}
This work was supported under the ISTI Rapid Response Research Call for "Hands-on Quantum Computing 2018" 
at Los Alamos National Laboratory.

Assigned: LA-UR-19-23387. Los Alamos National Laboratory is operated by Triad National Security, LLC for the National Nuclear Security Administration of U.S. Department of Energy under contract 89233218CNA000001.

\bibliography{references}
\bibliographystyle{abbrvnat}

\end{document}